\documentclass[12pt]{article}
\usepackage{amsfonts}
\begin{document}
\newsymbol\rtimes 226F
\newsymbol\ltimes 226E
\newcommand{\text}[1]{\mbox{{\rm #1}}}
\newcommand{\G}{{\bf G}}
\newcommand{\Rep}{\text{Rep}}
\newcommand{\gr}{\text{gr}}
\newcommand{\Fun}{\text{Fun}}
\newcommand{\Hom}{\text{Hom}}
\newcommand{\End}{\text{End}}
\newcommand{\FPdim}{\text{FPdim}}
\newcommand{\GL}{\text{GL}}
\newcommand{\Sp}{\text{Sp}}
\newcommand{\Ps}{\text{Ps}}
\newcommand{\Ad}{\text{Ad}}
\newcommand{\ASp}{\text{ASp}}
\newcommand{\APs}{\text{APs}}
\newcommand{\Rad}{\text{Rad}}
\newcommand{\Corad}{\text{Corad}}
\newcommand{\SuperVect}{\text{SuperVect}}
\newcommand{\Vect}{\text{Vect}}
\newcommand{\Spec}{\text{Spec}}
\newcommand{\tr}{\text{tr}}
\newcommand{\cH}{{\cal A}}
\newcommand{\cR}{{\cal R}}
\newcommand{\cJ}{\cal J}
\newcommand{\gd}{\delta}
\newcommand{\lan}{\langle}
\newcommand{\ran}{\rangle}
\newcommand{\itms}[1]{\item[[#1]]}
\newcommand{\nin}{\in\!\!\!\!\!/}
\newcommand{\g}{{\bf g}}
\newcommand{\sub}{\subset}
\newcommand{\cntd}{\subseteq}
\newcommand{\go}{\omega}
\newcommand{\Pa}{P_{a^\nu,1}(U)}
\newcommand{\fx}{f(x)}
\newcommand{\fy}{f(y)}
\newcommand{\gD}{\Delta}
\newcommand{\gl}{\lambda}
\newcommand{\gL}{\Lambda}
\newcommand{\half}{\frac{1}{2}}
\newcommand{\sto}[1]{#1^{(1)}}
\newcommand{\stt}[1]{#1^{(2)}}
\newcommand{\Z}{\hbox{\sf Z\kern-0.720em\hbox{ Z}}}
\newcommand{\singcolb}[2]{\left(\begin{array}{c}#1\\#2
\end{array}\right)}
\newcommand{\ga}{\alpha}
\newcommand{\gb}{\beta}
\newcommand{\gga}{\gamma}
\newcommand{\ul}{\underline}
\newcommand{\ol}{\overline}
\newcommand{\qed}{\kern 5pt\vrule height8pt width6.5pt depth2pt}
\newcommand{\Lrraro}{\Longrightarrow}
\newcommand{\Nb}{|\!\!/}
\newcommand{\NN}{{\rm I\!N}}
\newcommand{\bsl}{\backslash}
\newcommand{\gt}{\theta}
\newcommand{\op}{\oplus}
\newcommand{\C}{{\bf C}}
\newcommand{\Q}{{\bf Q}}
\newcommand{\Op}{\bigoplus}
\newcommand{\CR}{{\cal R}}
\newcommand{\grr}{\omega_1}
\newcommand{\ben}{\begin{enumerate}}
\newcommand{\een}{\end{enumerate}}
\newcommand{\ndiv}{\not\mid}
\newcommand{\bab}{\bowtie}
\newcommand{\hal}{\leftharpoonup}
\newcommand{\har}{\rightharpoonup}
\newcommand{\ot}{\otimes}
\newcommand{\OT}{\bigotimes}
\newcommand{\bwe}{\bigwedge}
\newcommand{\eps}{\varepsilon}
\newcommand{\gs}{\sigma}
\newcommand{\rbraces}[1]{\left( #1 \right)}
\newcommand{\bbox}{$\;\;\rule{2mm}{2mm}$}
\newcommand{\sbraces}[1]{\left[ #1 \right]}
\newcommand{\bbraces}[1]{\left\{ #1 \right\}}
\newcommand{\OO}{_{(1)}}
\newcommand{\TT}{_{(2)}}
\newcommand{\FF}{_{(3)}}
\newcommand{\minus}{^{-1}}
\newcommand{\CV}{\cal V}
\newcommand{\CVs}{\cal{V}_s}
\newcommand{\un}{U_q(sl_n)'}
\newcommand{\on}{O_q(SL_n)'}
\newcommand{\slq}{U_q(sl_2)}
\newcommand{\olq}{O_q(SL_2)}
\newcommand{\UU}{U_{(N,\nu,\go)}}
\newcommand{\HH}{{\mathcal H}}
\newcommand{\ct}{\centerline}
\newcommand{\bs}{\bigskip}
\newcommand{\qua}{\rm quasitriangular}
\newcommand{\ms}{\medskip}
\newcommand{\noin}{\noindent}
\newcommand{\mat}[1]{$\;{#1}\;$}
\newcommand{\raro}{\rightarrow}
\newcommand{\map}[3]{{#1}\::\:{#2}\raro{#3}}
\newcommand{\alg}{{\rm Alg}}
\def\newtheorems{\newtheorem{theorem}{Theorem}[subsection]
                 \newtheorem{corollary}[theorem]{Corollary}
                 \newtheorem{proposition}[theorem]{Proposition}
                 \newtheorem{lemma}[theorem]{Lemma}
                 \newtheorem{definition}[theorem]{Definition}
                 \newtheorem{Theorem}{Theorem}[section]
                 \newtheorem{Corollary}[Theorem]{Corollary}
                 \newtheorem{Proposition}[Theorem]{Proposition}
                 \newtheorem{Lemma}[Theorem]{Lemma}
                 \newtheorem{Definition}[Theorem]{Definition}
                 \newtheorem{Example}[Theorem]{Example}
                 \newtheorem{Remark}[Theorem]{Remark}
                 \newtheorem{claim}[theorem]{Claim}
                 \newtheorem{sublemma}[theorem]{Sublemma}
                 \newtheorem{example}[theorem]{Example}
                 \newtheorem{remark}[theorem]{Remark}
                 \newtheorem{question}[theorem]{Question}
                 \newtheorem{Question}[Theorem]{Question}
                 \newtheorem{Conjecture}[Theorem]{Conjecture}
                 \newtheorem{ttheorem}{Theorem}[subsubsection]
                 \newtheorem{qquestion}[ttheorem]{Question}
                 \newtheorem{cconjecture}[ttheorem]{Conjecture}
                 \newtheorem{eexample}[ttheorem]{Example}
                 \newtheorem{rremark}[ttheorem]{Remark}
                 \newtheorem{ccorollary}[ttheorem]{Corollary}
                 \newtheorem{pproposition}[ttheorem]{Proposition}
                 \newtheorem{llemma}[ttheorem]{Lemma}
                 \newtheorem{ddefinition}[ttheorem]{Definition}
}
\newtheorems
\newcommand{\proof}{\par\noindent{\bf Proof:}\quad}
\newcommand{\dmatr}[2]{\left(\begin{array}{c}{#1}\\
                            {#2}\end{array}\right)}
\newcommand{\doubcolb}[4]{\left(\begin{array}{cc}#1&#2\\
#3&#4\end{array}\right)}
\newcommand{\qmatrl}[4]{\left(\begin{array}{ll}{#1}&{#2}\\
                            {#3}&{#4}\end{array}\right)}
\newcommand{\qmatrc}[4]{\left(\begin{array}{cc}{#1}&{#2}\\
                            {#3}&{#4}\end{array}\right)}
\newcommand{\qmatrr}[4]{\left(\begin{array}{rr}{#1}&{#2}\\
                            {#3}&{#4}\end{array}\right)}
\newcommand{\smatr}[2]{\left(\begin{array}{c}{#1}\\
                            \vdots\\{#2}\end{array}\right)}

\newcommand{\ddet}[2]{\left[\begin{array}{c}{#1}\\
                           {#2}\end{array}\right]}
\newcommand{\qdetl}[4]{\left[\begin{array}{ll}{#1}&{#2}\\
                           {#3}&{#4}\end{array}\right]}
\newcommand{\qdetc}[4]{\left[\begin{array}{cc}{#1}&{#2}\\
                           {#3}&{#4}\end{array}\right]}
\newcommand{\qdetr}[4]{\left[\begin{array}{rr}{#1}&{#2}\\
                           {#3}&{#4}\end{array}\right]}

\newcommand{\qbracl}[4]{\left\{\begin{array}{ll}{#1}&{#2}\\
                           {#3}&{#4}\end{array}\right.}
\newcommand{\qbracr}[4]{\left.\begin{array}{ll}{#1}&{#2}\\
                           {#3}&{#4}\end{array}\right\}}

\title{{\bf On the Classification of Finite-Dimensional
Triangular Hopf Algebras}}
\author{Shlomo Gelaki\\Technion-Israel
Institute of Technology\\
Department of Mathematics\\
Haifa 32000, Israel\\
email: gelaki@math.technion.ac.il}
\maketitle

\section{Introduction}
A fundamental problem in the theory of
Hopf algebras is
the classification and construction of finite-dimensional
quasitriangular Hopf algebras $(A,R)$ over an
algebraically closed field $k.$ Quasitriangular Hopf
algebras constitute a very important class of Hopf
algebras, which were introduced by Drinfeld
[Dr1] in order to supply solutions to the
quantum Yang-Baxter equation that arises in
mathematical physics. Quasitriangular Hopf algebras are the 
Hopf algebras whose finite-dimensional representations form a 
braided rigid tensor category, which naturally relates
them to low dimensional topology. Furthermore, Drinfeld
showed that {\em any} finite-dimensional Hopf algebra can be
embedded in a finite-dimensional quasitriangular Hopf
algebra, known now as its Drinfeld double or quantum
double.
However, this intriguing problem turns out to be extremely
hard and it is still widely open. One can hope that 
resolving this problem first in the {\em triangular} case
would contribute to the understanding of the general
problem. 

Triangular Hopf algebras are the Hopf algebras whose
representations
form a symmetric tensor category. In that sense, they are
the class of Hopf algebras closest to group algebras.
The structure of triangular Hopf algebras is far from trivial,
and yet is more tractable than that of general Hopf
algebras, due to their proximity to groups and Lie algebras.
This makes triangular Hopf algebras an excellent testing
ground for
general Hopf algebraic ideas, methods and conjectures.
A general classification
of triangular Hopf algebras is not known yet.
However, there are two classes
that are relatively well understood.
One of them is semisimple triangular Hopf algebras over $k$
(and cosemisimple if the characteristic of $k$ is positive)
for which a complete classification is given in [EG1,EG4].
The key theorem about such Hopf algebras states that each
of them is obtained by twisting a group algebra of a finite group
[EG1, Theorem 2.1] (see also [G5]). 

Another important class of Hopf algebras is that of
{\em pointed} ones. These are Hopf algebras whose all
simple comodules are $1-$dimensional.
Theorem 5.1 in [G4] (together with [AEG, Theorem 6.1]) gives a
classification of
{\em minimal} triangular pointed Hopf algebras.

Recall that a finite-dimensional algebra is called {\em
basic} if all of its simple modules are $1-$dimensional
(i.e. if its dual is a pointed coalgebra).
The same Theorem 5.1 of [G4] gives a classification
of minimal triangular basic Hopf algebras, since the dual of a
minimal triangular Hopf algebra is again minimal triangular.

Basic and semisimple Hopf algebras share a common property. Namely, 
the Jacobson radical $\Rad(A)$ of such a Hopf algebra $A$ is a
Hopf ideal,
and hence the quotient $A/\Rad(A)$ (the semisimple part)
is itself a Hopf algebra. The representation-theoretic
formulation
of this property is: The tensor product of two simple
$A$-modules
is semisimple. A remarkable classical theorem of Chevalley
[C, p.88] states that, in characteristic $0,$ this property holds for
the group algebra of any (not necessarily finite) group. So we 
called this property of $A$ {\bf the Chevalley property} [AEG].

In [AEG] it was proved that any finite-dimensional triangular Hopf
algebra with the Chevalley property
is obtained by twisting a finite-dimensional triangular
Hopf algebra with $R-$matrix of rank $\le 2,$
and that any finite-dimensional triangular Hopf algebra with
$R-$matrix of rank $\le 2$
is a suitable modification of a finite-dimensional 
cocommutative Hopf   
superalgebra (i.e. the group algebra of a finite supergroup).
On the other hand, by a
theorem of Kostant [Ko], a finite supergroup is a semidirect
product of a finite group with an odd vector space on which
this group acts. Moreover, the converse result
that any such Hopf algebra does have the Chevalley property is
also proved in [AEG]. As a corollary, we proved that any
finite-dimensional
triangular Hopf algebra whose coradical is a Hopf subalgebra
(e.g. pointed) is obtained by twisting a triangular
Hopf algebra with $R-$matrix of rank $\le 2$.

The purpose of this paper is to present all that is currently
known to us about the classification and construction of
finite-dimensional triangular Hopf algebras, and to explain
the mathematics and ideas involved in this theory.

The paper is organized as follows. In Section 2 we review
some necessary material from the theory of Hopf algebras.
In particular the important notion of a twist for Hopf
algebras, which was introduced by Drinfeld [Dr1].

In Section 3 we explain in details the theory of Movshev
on twisting in group algebras of finite groups [Mov].
The results of [EG4,EG5] (described in Sections 4 and 5
below) rely, among other things, on this theory in an 
essential way. 

In Section 4 we concentrate on the theory of triangular
semisimple and cosemisimple Hopf algebras. We
first describe the classification and construction of
triangular semisimple and cosemisimple Hopf algebras
over {\em any}
algebraically closed field $k,$ and then describe some 
of the consequences of the classification theorem, in
particular the one concerning the existence of
grouplike elements in triangular semisimple and
cosemisimple Hopf algebras over $k$ [EG4].
The classification uses, among other things, 
Deligne's theorem on Tannakian
categories [De1] in
an essential way. We refer the reader to [G5] for
a detailed discussion of this aspect.
The proof of the existence of
grouplike elements relies on a theorem from [HI] on
central type groups being solvable, which is
proved using the classification of finite simple groups.
The classification in positive characteristic relies also
on the lifting functor from [EG5].

In Section 5 we concentrate on the dual
objects of Section 4; namely, on semisimple and
cosemisimple {\em cotriangular} Hopf algebras over $k,$
studied in [EG3]. We describe the
representation theory of such Hopf algebras, and in
particular obtain that Kaplansky's 6th conjecture
[Kap] holds for them (i.e. they are of Frobenius type).

In Section 6 we concentrate on the pointed case, studied in
[G4] and [AEG, Theorem 6.1]. The main result in this case is
the classification of minimal triangular pointed Hopf
algebras. 

In Section 7 we generalize and concentrate on the
classification of finite-dimensional triangular Hopf algebras
with the
Chevalley property, given in [AEG]. We note that
similarly to the
case of semisimple Hopf algebras, the proof of the main result of [AEG]
is based on Deligne's theorem [De1]. In fact, we used Theorem 2.1 of
[EG1] to prove the main result of this paper.

In Section 8 we conclude the paper with a list of relevant
questions raised in [AEG] and [G4].

Throughout the paper the ground field
$k$ is assumed to be algebraically closed.
The symbol $\C$ will always denote the field of complex numbers.
For a Hopf (super)algebra $A,$ $\G(A)$ will denote its group
of grouplike elements.

\noin
{\bf Acknowledgment} The work described in
Sections 3-5 is joint with Pavel
Etingof, whom I am grateful to for his
help in reading the manuscript. The work described in
Subsection 4.5 is joint also with Robert Guralnick and 
Jan Saxl. The work described in
Section 7 is joint with Nicholas Andruskiewitsch and Pavel 
Etingof. 

\section{Preliminaries}
In this Section we recall the necessary background
needed for
this paper. We refer the reader to the books [ES,Kass,Mon,Sw] 
for the general theory of Hopf algebras and quantum
groups. 

\subsection{Quasitriangular Hopf algebras}
We recall Drinfeld's notion of a (quasi)triangular Hopf
algebra [Dr1].
Let $(A,m,1,\gD,\varepsilon,S)$ be a
finite-dimensional Hopf
algebra over $k,$ and let
$R=\sum_i a_i\ot b_i\in A\ot A$ be an invertible
element. Define a
linear map $f_R:A^*\raro A$ by
$f_R(p)=\sum_i <p,a_i>b_i$ for
$p\in A^*$. 
The tuple $(A,m,1,\gD,\varepsilon,S,R)$ is said to be a 
{\em
quasitriangular} Hopf algebra if the  
following axioms hold:
\begin{equation}\label{QT.1}
(\gD \ot Id)(R)=R_{13}R_{23},\, (Id\ot \gD)(R)=R_{13}R_{12}
\end{equation}
where $Id$ is the identity map of $A,$ and
\begin{equation}\label{QT.5}
\gD^{cop}(a)R=R\gD(a)\;\text{for any}\;a\in A;   
\end{equation}

\noin 
or equivalently, if $f_R:A^*\raro A^{cop}$ is a Hopf 
algebra map and (\ref{QT.5}) is satisfied. The element $R$ is 
called an $R-$matrix. Observe that
using Sweedler's notation for the comultiplication [Sw],
(\ref{QT.5}) is equivalent to
\begin{equation}\label{qt5}   
\sum <p_{(1)},a_{(2)}>a_{(1)}f_{_{R}}(p_{(2)})=\sum
<p_{(2)},a_{(1)}> f_{_{R}}(p_{(1)})a_{(2)}
\end{equation}
for any $p\in A^*$ and $a\in A.$

A quasitriangular Hopf algebra $(A,R)$ is called {\em
triangular} if $R^{-1}=R_{21};$ or equivalently, 
if $f_R*f_{R_{21}}=\varepsilon$ in
the convolution algebra $\text{Hom}_k(A^*,A),$ i.e.
(using Sweedler's notation again)
\begin{equation}\label{mapt}
\sum f_R(p_{(1)})f_{R_{21}}(p_{(2)})=<p,1>1\;\text{for
any}\;p\in A^*.
\end{equation}

Let 
\begin{equation}\label{drinelm}
u:=\sum_i S(b_i)a_i 
\end{equation}
be the {\em Drinfeld element} of $(A,R).$ Drinfeld
showed
[Dr2] that $u$ is invertible and that
\begin{equation}\label{ssqu}
S^2(a)=uau^{-1}\;\text{for any}\;a\in A.
\end{equation}
He also showed that $(A,R)$ is triangular if and only
if $u$ is a grouplike element [Dr2].

Suppose further that $(A,m,1,\Delta,\varepsilon,S,R)$
is {\em semisimple and cosemisimple} over $k.$ 
\begin{lemma}\label{u1}
The Drinfeld element $u$ is central, and 
\begin{equation}\label{usu}
u=S(u).
\end{equation}
\end{lemma}
\proof
By [LR1] in characteristic
$0,$ and by [EG5, Theorem 3.1] in positive
characteristic, $S^2=I.$ Hence by (\ref{ssqu}), $u$ is
central.
Now, we have $(S\ot S)(R)=R$ [Dr2], so $S(u)=
\sum_i S(a_i)S^2(b_i)=\sum_i a_iS(b_i)$. This shows that
$\tr(u)=\tr(S(u))$ in
every irreducible representation of $A$. But $u$ and
$S(u)$ are central,
so they act as scalars in this representation, which proves
(\ref{usu}). \qed
\begin{lemma}\label{u2}
In particular, 
\begin{equation}\label{usq}
u^2=1.
\end{equation}
\end{lemma}
\proof
Since $S(u)=u^{-1},$ the result
follows from (\ref{usu}). \qed

Let us demonstrate that it is always possible 
to replace $R$ with a new $R-$matrix $\tilde R$ so
that the new Drinfeld element $\tilde u$ equals $1.$
Indeed, if $k$ does not have characteristic $2,$ set
\begin{equation}\label{ru}
R_u:=\frac{1}{2}(1\ot 1+1\ot u+u\ot 1-u\ot u).
\end{equation}
If $k$ is of characteristic $2$ (in which case $u=1$
by semisimplicity), set $R_u:=1.$ Set $\tilde
R:=RR_u.$
\begin{lemma}\label{u3}
$(A,\tilde R)$ is a triangular semisimple and
cosemisimple Hopf algebra with Drinfeld 
element $1.$ 
\end{lemma}
\proof Straightforward. \qed

This observation allows to reduce questions about
triangular semisimple and cosemisimple
Hopf algebras over $k$ to the case when the Drinfeld
element is $1.$

Let $(A,R)$ be {\em any} triangular Hopf algebra over $k.$
Write $R=\sum_{i=1}^na_i\ot b_i$ in the shortest possible way, 
and let $A_m$ be the Hopf subalgebra of $A$ generated by the
$a_i$'s and $b_i$'s. Following [R2], we will call $A_m$ the 
{\em minimal part} of $A.$ We will call $n=\dim(A_m)$ the {\em
rank} of the
$R-$matrix $R.$
It is straightforward to
verify that
the corresponding map $f_R:A_m^{*cop}\raro A_m$ defined by 
$f_R(p)=(p\ot I)(R)$ is a Hopf algebra isomorphism. This 
property of minimal triangular Hopf algebras will
play a central
role in our study of the pointed case (see Section 6 
below). It
implies in particular that
$\G(A_m)\cong \G((A_m)^*),$ and hence that the group
$\G(A_m)$ is
{\em abelian} (see e.g. [G2]). Thus, $\G(A_m)\cong
\G(A_m)^\vee $ (where
$\G(A_m)^\vee $ denotes the character group of
$\G(A_m)$), and we can identify the Hopf algebras
$k[\G(A_m)^\vee ]$ and $k[\G(A_m)]^*.$  
Also, if $(A,R)$ is minimal triangular and pointed
then $f_R$
being an isomorphism implies that $A^*$ is pointed as well.

Note that if $(A,R)$ is (quasi)triangular and
$\pi:A\raro A'$ is a surjective map of Hopf algebras, then
$(A',R')$ is
(quasi)triangular as well, where $R':=(\pi\ot \pi)(R).$ 

\subsection{Hopf superalgebras}

\subsubsection{Supervector spaces}

We start by recalling the definition of the category of
supervector spaces.
A Hopf algebraic way to define this category is as follows.
Let us assume that $k=\C.$

Let $u$ be the generator of the group $\Z_2$ of two
elements, and let $R_u\in \C[\Z_2]\otimes \C[\Z_2]$ be as 
in (\ref{ru}). Then
$(\C[\Z_2],R_u)$ is a
minimal triangular Hopf algebra.

\begin{ddefinition} The category of supervector spaces
over $\C$
is the symmetric tensor category $\Rep(\C[\Z_2],R_u)$
of representations of the triangular Hopf algebra
$(\C[\Z_2],R_u).$ This category will be denoted by
$\SuperVect.$
\end{ddefinition}

For $V\in \SuperVect$ and $v\in V$, we say that $v$ is
even
if $uv=v$ and odd if $uv=-v$. The set of even vectors in
$V$ is
denoted by $V_0$ and the set of odd vectors by $V_1$, so
$V=V_0\oplus V_1$.
We define the parity of a vector $v$ to be $p(v)=0$ if $v$ is
even and $p(v)=1$ if $v$ is odd (if $v$ is neither odd
nor even,
$p(v)$ is not defined).

Thus, as an ordinary tensor category, SuperVect is equivalent
to the category of representations of $\Z_2$, but the
commutativity
constraint is different
from that of $\Rep(\Z_2)$ and equals $\beta:=R_uP$, where
$P$ is the
permutation of components.
In other words, we have
\begin{equation}\label{symm}
\beta (v\otimes w)=(-1)^{p(v)p(w)}w\otimes v,
\end{equation}
where both $v,w$ are either even or odd.

\subsubsection{Hopf superalgebras}

Recall that in any symmetric (more generally,
braided) tensor category,
one can define an algebra, coalgebra,
bialgebra, Hopf algebra, triangular Hopf algebra, etc,
to be an object of this category equipped with
the usual structure maps
(morphisms in this category),
subject to the same axioms as in the usual case.
In particular, any of these algebraic structures
in the category SuperVect is usually identified by the prefix
``super''. For example:

\begin{ddefinition}
A Hopf superalgebra is a Hopf algebra
in
$\SuperVect.$
\end{ddefinition}

More specifically, a Hopf superalgebra $\cH$
is an ordinary
$\Z_2$-graded associative unital algebra with
multiplication $m$,
equipped with a coassociative map
$\Delta:\cH\raro
\cH\ot \cH$ (a morphism in $\SuperVect$) which is
multiplicative in
the super-sense, and with a
counit and antipode satisfying the standard axioms.
Here multiplicativity in the super-sense
means that $\Delta$ satisfies the relation
\begin{equation}\label{sd}
\Delta(ab)=\sum (-1)^{p(a_2)p(b_1)}a_1b_1\ot a_2b_2
\end{equation}
for all $a,b\in \cH$
(where $\Delta(a)=\sum a_1\ot a_2$, $\Delta(b)=\sum
b_1\ot b_2$).
This is because the tensor product of
two algebras $A,B$ in $\SuperVect$ is defined to be $A\ot
B$ as a vector space, with multiplication
\begin{equation}\label{sm}
(a\ot b)(a'\ot b'):=(-1)^{p(a')p(b)}a a'\ot
bb'.
\end{equation}

\begin{rremark} {\rm Hopf superalgebras appear in
[Ko], under the name of ``graded Hopf algebras''.
}
\end{rremark}

Similarly, a
(quasi)triangular Hopf superalgebra $(\cH,\cR)$ is a Hopf
superalgebra
with an $R-$matrix (an {\em even} element $\cR\in \cH\otimes
\cH$)
satisfying the usual axioms. As in the even case, an
important role is played by the Drinfeld element $u$ of
$(\cH,\cR)$:
\begin{equation}\label{du}
u:=m\circ \beta \circ (Id\ot S)(\cR).
\end{equation}
For instance, $(\cH,\cR)$ is triangular if and only if $u$
is a grouplike element of $\cH.$

As in the even case, the tensorands of the $R-$matrix
of a (quasi)triangular Hopf superalgebra $\cH$ generate a
finite-dimensional
sub Hopf superalgebra ${\cH}_m$, called the {\em minimal
part of} $\cH$ (the proof does not differ essentially
from the proof of the analogous fact for Hopf algebras).
A (quasi)triangular Hopf superalgebra is said to be minimal if
it coincides with its minimal part. The dimension of the minimal part
is the {\it rank} of the $R-$matrix.

\subsubsection{Cocommutative Hopf superalgebras}

\begin{ddefinition}
We will say that a Hopf superalgebra $\cH$ is
commutative (resp. cocommutative) if $m=m\circ
\beta$
(resp. $\Delta=\beta\circ \Delta$).
\end{ddefinition}

\begin{eexample}\label{supergr}
{\rm {\bf [Ko]} Let $G$ be a group, and $\g$ a Lie
superalgebra with an action of
$G$ by automorphisms of Lie superalgebras.
Let ${\cal A}:=\C[G] \ltimes \text{U}(\g),$ where
$\text{U}(\g)$
denotes the universal
enveloping algebra of $\g$.
Then $\cH$ is a
cocommutative Hopf
superalgebra, with
$\Delta(x)=x\ot 1+1\ot x,$ $x\in \g$, and $\Delta(g)=g\ot g,$
$g\in G$. In this Hopf superalgebra, we have
$S(g)=g^{-1}$, $S(x)=-x$, and in particular $S^2=Id$.

The Hopf superalgebra $\cH$ is finite-dimensional if and
only if
$G$ is finite, and $\g$ is finite-dimensional and purely
odd (and
hence commutative).
Then
${\cal A}=\C[G]\ltimes \Lambda V$, where $V=\g$ is an odd
vector space
with a $G$-action. In this case, $\cH^*$ is a commutative
Hopf
superalgebra.
}
\end{eexample}

\begin{rremark} {\rm We note that as in the even case, it
is convenient to think about $\cH$ and $\cH^*$ in geometric terms.
Consider, for instance, the finite-dimensional case. In this case,
it is useful to think of the ``affine algebraic supergroup''
$\tilde G:=G\ltimes V$. Then one can regard $\cH$ as the group
algebra $\C[\tilde G]$ of this supergroup, and $\cH^*$ as its
function algebra $F(\tilde G)$. Having this in mind, we will call
the algebra $\cH$ {\bf a supergroup algebra}. }
\end{rremark}

It turns out that like in the even case,
any cocommutative Hopf superalgebra is of the type
described
in Example \ref{supergr}. Namely, we have the following
theorem.

\begin{ttheorem}\label{kostant} {\bf ([Ko], Theorem 3.3)}
Let $\cH$ be a cocommutative Hopf superalgebra over $\C.$
Then $\cH=\C[{\bf G}(\cH)]\ltimes
\text{U}(P(\cH)),$
where $\text{U}(P(\cH))$ is the universal enveloping
algebra
of the
Lie superalgebra of primitive elements of $\cH,$
and ${\bf G}(\cH)$ is the group of grouplike elements of
$\cH.$
\end{ttheorem}

In particular, in the finite-dimensional case we get:

\begin{ccorollary}\label{kostantf}
Let $\cH$ be a finite-dimensional cocommutative
Hopf superalgebra over $\C.$
Then ${\cal A}=\C[{\bf G}({\cal A})]\ltimes \Lambda V,$
where $V$ is the space of primitive elements of $\cH$
(regarded as an odd vector space) and ${\bf G}({\cal A})$ is the
finite group of grouplikes of ${\cal A}.$
In other words, $\cH$ is a supergroup algebra.
\end{ccorollary}

\subsection{Twists}
Let $(A,m,1,\gD,\varepsilon,S)$ be a Hopf algebra over a
field $k.$ We recall Drinfeld's notion of a {\em twist}
for $A$ [Dr1]. 

\begin{definition}\label{t}
A quasitwist for $A$ is an
element $J\in A\ot A$ which satisfies
\begin{equation}\label{t1}
(\Delta\ot Id)(J)(J\ot 1)=(Id\ot \Delta)(J)(1\ot J)\;\;
and \;\;
(\varepsilon\ot
Id)(J)=(Id\ot \varepsilon)(J)=1.
\end{equation}
An invertible quasitwist for $A$ is called a twist.
\end{definition}

Given a twist $J$ for $A$, one can define a new Hopf algebra
structure $(A^J,m,1,\Delta^J,\varepsilon,S^J)$ on the
algebra $(A,m,1)$ 
as follows. The coproduct is determined by
\begin{equation}\label{t2}
\Delta^J(a)=J^{-1}\Delta(a)J\;\text{for any}\;a\in A,
\end{equation}
and the antipode is determined by
\begin{equation}\label{tant}
S^J(a)=Q^{-1}S(a)Q\;\text{for any}\;a\in A,
\end{equation}
where $Q:=m\circ(S\ot Id)(J).$ If $A$ is (quasi)triangular
with the
universal $R-$matrix $R,$ then so is $A^J,$ with
the universal $R-$matrix $R^J:=J_{21}^{-1}RJ.$ 

\begin{example}\label{abt}
{\rm Let $G$ be a finite {\em abelian} group, and $G^\vee $ its
character group. Then the set of twists for
$A:=k[G]$ is in one to one correspondence with the set
of $2-$cocycles $c$ of $G^\vee $ with coefficients in $k^*,$
such that $c(0,0)=1.$
Indeed, let $J$ be a twist for $A,$ and define
$c:G^\vee \times G^\vee \raro k^*$ via
$c(\chi,\psi):=(\chi\ot\psi)(J).$ Then it is
straightforward to verify that $c$ is a $2-$cocycle
of $G^\vee $ (see e.g. [Mov, Proposition 3]), and that
$c(0,0)=1.$

Conversely, let $c:G^\vee \times G^\vee \raro k^*$ be a
$2-$cocycle of
$G^\vee $ with coefficients in $k^*,$ such that $c(0,0)=1.$
Note that the $2-$cocycle condition implies that $c(0,\chi)=1=
c(\chi,0)$ for all $\chi\in G^\vee.$ For $\chi\in G^\vee,$
let $E_{\chi}:=|G|^{-1}\sum_{g\in G}\chi(g)g$ be the
associated idempotent of $A.$ Then
it is straightforward to verify that
$J:=\sum_{\chi,\psi\in G^\vee }c(\chi,\psi)E_{\chi}\ot
E_{\psi}$ is a twist for $A$ (see e.g. [Mov,
Proposition 3]). Moreover it is easy to check that the
above two assignments are inverse to each other.
}
\end{example}

\begin{remark}\label{nonab}
{\rm Unlike for finite abelian groups, 
the study of twists for finite non-abelian groups is
much more involved. This was done in [EG2,EG4,Mov] (see
Section 4 below).}
\end{remark}

If $J$ is a (quasi)twist for $A$ and $x$ is an invertible
element of $A$ such that $\varepsilon(x)=1, $ then
\begin{equation}\label{jgauge}
J^x:=\Delta(x)J(x^{-1}\otimes x^{-1})
\end{equation}
is also a (quasi)twist for $A.$
We will call the (quasi)twists $J$ and $J^x$ {\em gauge
equivalent}.
Observe that if $(A,R)$ is a (quasi)triangular Hopf algebra, then
the map $(A^J,R^J)\raro (A^{J^x},R^{J^x})$ determined
by $a\mapsto xax^{-1}$ is
an isomorphism of (quasi)triangular Hopf algebras.

Let $A$ be a group algebra of a finite group. We will say 
that a twist
$J$ for $A$ is {\em minimal} if the
right (and left) components of
the $R-$matrix $R^J:=J_{21}^{-1}J$
span $A,$ i.e. if the corresponding triangular Hopf algebra
$(A^J,J_{21}^{-1}J)$ is minimal.

A twist for a Hopf algebra in {\em any symmetric tensor category}
is defined in the same way as in the usual case. For instance, if $\cH$ 
is a Hopf superalgebra then a twist for $\cH$ is an invertible {\em 
even} element ${\cal J}\in \cH\ot \cH$ satisfying (\ref{t1}).

\subsection{Projective representations and central extensions}

Here we recall some basic facts about projective
representations
and central extensions. They can be found in textbooks, e.g.
[CR, Section 11E].

A projective representation over $k$ of a group $H$ is a
vector space $V$
together with a homomorphism of groups $\pi_{_V}:H\raro
\text{PGL}(V),$ where
$\text{PGL}(V)\cong \text{GL}(V)/k^*$ is the projective linear
group.

A linearization of a projective representation $V$ of $H$
is a central extension $\hat H$ of $H$ by a central subgroup
$\zeta$ together with a linear representation
$\tilde\pi_{_V}:\hat H\to \text{GL}(V)$ which descends to
$\pi_{_V}$.
If $V$ is a finite-dimensional projective representation
of $H$
then there exists a linearization of $V$ such that
$\zeta$ is finite
(in fact, one can make $\zeta=\Z/(\text{dim}(V))\Z$).

Any projective representation $V$ of $H$ canonically defines a
cohomology class
$[V]\in H^2(H,k^*)$. The representation $V$ can be lifted to
a linear representation of $H$ if and only if $[V]=0$.

\subsection{Pointed Hopf algebras} 
The Hopf algebras which are studied in Section 6 are
pointed. 
Recall that a Hopf algebra $A$ is {\em pointed} if its simple
subcoalgebras are 
all $1-$dimensional or equivalently (when $A$ is
finite-dimensional) if 
the irreducible representations of $A^*$ are all
$1-$dimensional (i.e. $A^*$ is basic). 
For any $g,h\in \G(A),$ we denote the vector space of
$g:h$ {\em skew primitives}
of $A$ by $P_{g,h}(A):=\{x\in A|\gD(x)=x\ot g+h\ot x\}.$
Thus the classical {\em primitive} elements of $A$ are
$P(A):=P_{1,1}(A).$
The element $g-h$ is always $g:h$ skew primitive. Let
$P'_{g,h}(A)$ denote a complement of $sp_k\{g-h\}$ in  
$P_{g,h}(A)$. Taft-Wilson theorem [TW] states that  
the first term $A_1$ of the coradical filtration of $A$ is
given by:
\begin{equation}\label{wilson}
A_1=k[\G(A)]\bigoplus \left(\Op_{g,h\in \G(A)}
P'_{g,h}(A)\right).
\end{equation}
In particular, if $A$ is {\em not} cosemisimple then there
exists $g\in \G(A)$ such that $P'_{1,g}(A)\ne 0.$ 

If $A$ is a Hopf algebra over the field $k,$ which is
generated
as an algebra by a subset $S$ of $\G(A)$ and by $g:g'$ skew
primitive
elements, where $g,g'$ run over $S,$ then $A$ is pointed
and $\G(A)$ is
generated as a group by $S$ (see e.g. [R4, Lemma 1]). 

\section{Movshev's theory on the algebra associated
with a twist}
In this section we describe Movshev's theory on
twisting in group algebras of finite groups [Mov]. Our
classification theory of triangular semisimple and
cosemisimple 
Hopf algebras [EG4] (see Section 4 below), and our  
study of the representation theory of cotriangular 
semisimple and cosemisimple Hopf algebras [EG3] (see 
Section 
5 below) rely, among other things, 
on this theory in an essential way.

Let $k$ be an algebraically closed field whose
characteristic 
is relatively prime to $|G|.$
Let $A:=k[G]$ be the group algebra of a finite group
$G,$ equipped with
the usual multiplication,
unit, comultiplication, counit and antipode, denoted by
$m,$ $1,$ $\Delta,$ $\varepsilon$ and $S$ respectively. 
Let $J\in A\ot A.$ 
Movshev had the following nice idea of characterizing
quasitwists [Mov].
Let $(A_J,\Delta_J,\varepsilon)$ where $A_J=A$ as
vector spaces, and $\Delta_J$ is the map
\begin{equation}\label{mjr}
\Delta_J:A\raro A\ot A,\;a\mapsto \Delta(a)J.
\end{equation}
\begin{Proposition}\label{movelem}
$(A_J,\Delta_J,\varepsilon)$ is a coalgebra if and only if
$J$ is a quasitwist for $A.$
\end{Proposition}
\proof Straightforward. \qed

Regard $A$ as the left regular representation of
$G.$ Then $(A_J,\Delta_J,\varepsilon)$ is a 
$G-$coalgebra (i.e. $\Delta_J(ga)=(g\ot g)\Delta_J(a)$
and $\varepsilon(ga)=\varepsilon(a)$ for all $g\in G,$
$a\in A$). In fact, we have the following important result.
\begin{Proposition} {\bf [Mov, Proposition 
5]}\label{movprop5}
Suppose that $(C,\tilde \Delta,\tilde\varepsilon)$ is a
$G-$coalgebra which is isomorphic to the regular
representation of $G$ as a $G-$module. Then there exists
a quasitwist $J\in A\ot A$ such that 
$(C,\tilde \Delta,\tilde\varepsilon)$ and
$(A,\Delta_J,\varepsilon)$ are isomorphic as
$G-$coalgebras. Moreover, $J$ is unique up to gauge
equivalence.
\end{Proposition}
\proof
We can choose an element $\lambda\in C$ such
that the set $\{g\cdot \lambda|g\in G\}$ forms a basis
of
$C,$ and $\tilde\varepsilon(\lambda)=1.$ Now, write
$\tilde\Delta(\lambda)=\sum_{a,b\in G} \gamma(a,b) a\cdot
\lambda\ot b\cdot \lambda,$ and set 
\begin{equation}\label{regj}
J:=\sum_{a,b\in G} \gamma(a,b) a\ot b\in A\ot A.
\end{equation}
We have to show that $J$ is a quasitwist for $A.$
Indeed, let
$f:A\raro C$ be determined by $f(a)=a\cdot \lambda.$
Clearly, $f$ is an isomorphism of
$G-$modules which satisfies $\tilde\Delta(f(a))=(f\ot
f)\Delta_J(a),$ $a\in A.$ Therefore 
$(A_J,\Delta_J,\varepsilon)$ is a coalgebra, which is
equivalent to saying that $J$ is a quasitwist by
Proposition \ref{movelem}. This proves the first claim.

Suppose that $(A_{J'},\Delta_{J'},\varepsilon)$ and
$(A_J,\Delta_J,\varepsilon)$ are isomorphic as 
$G-$coalgebras via $\phi:A \raro A.$ We have to show 
that $J,$ $J'$ are gauge equivalent. Indeed, $\phi$ is
given by
right multiplication by an invertible element $x\in A,$
$\phi(a)=ax.$ On one hand, $(\phi\ot
\phi)(\Delta_J(1))=J(x\ot x),$ and on
the other hand, 
$\Delta_{J'}(\phi(1))=\Delta(x)J'.$
The equality between the two right hand sides implies the
desired result. \qed

We now focus on the dual algebra $(A_J)^*$ of the
coalgebra $(A_J,\Delta_J,\varepsilon),$ and summarize
Movshev's results about it [Mov]. Note that
$(A_J)^*$ is a $G-$algebra which is isomorphic to the
regular representation of $G$ as a $G-$module.
\begin{Proposition} {\bf [Mov, Propositions
6,7]}\label{movprop67}
The following hold:
\ben
\item The algebra $(A_J)^*$ is semisimple. 
\item There exists a subgroup $St$ of $G$ (the
stabilizer of 
a maximal two sided ideal $I$ of $(A_J)^*$) such that
$(A_J)^*$ is isomorphic to the algebra of functions from
the set $G/St$ to the matrix algebra
$M_{|St|^{1/2}}(k).$
\een
\end{Proposition}

Note that, in particular, the group $St$ acts on the
matrix algebra $(A_J)^*/I\cong
M_{|St|^{1/2}}(k).$ Hence this algebra defines a projective
representation $T:St\raro \text{PGL}(|St|^{1/2},k)$
(since $\text{Aut}(M_{|St|^{1/2}}(k))=\text{PGL}(|St|^{1/2},k)$).
\begin{Proposition} {\bf [Mov, Propositions
8,9]}\label{movprop89}
$T$ is irreducible, and the associated $2-$cocycle $c:
St\times St\raro k^*$ is nontrivial.
\end{Proposition}

Consider now the twisted group algebra $k[St]^c.$
This algebra has a basis $\{X_g|g\in St\}$ with relations
$X_g X_h=c(g,h)X_{gh},$ and a natural
structure as a $St-$algebra given by \linebreak
$a\cdot X_g:=X_a X_g (X_a)^{-1}$ for all $a\in St$
(see also [Mov, Proposition 10]). Recall that $c$ is
called {\em nondegenerate} if for all $1\ne g\in St,$
the map $C_{St}(g)\raro k^*,$ $m\mapsto c(m,g)/c(g,m)$
is a nontrivial homomorphism of the centralizer of $g$
in $St$ to $k^*.$ In [Mov, Propositions 11,12]
Movshev reproduces the following well known criterion
for $k[St]^c$ to be a simple algebra (i.e. isomorphic
to the matrix algebra $M_{|St|^{1/2}}(k)$). 
\begin{Proposition}\label{movprop1112}
The twisted group algebra $k[St]^c$ is simple if and only
if $c$ is nondegenerate. Furthermore, if this is the
case, then
$k[St]^c$ is isomorphic to the regular representation
of $St$ as a $St-$module.
\end{Proposition}

Assume $c$ is nondegenerate. By Proposition
\ref{movprop5}, the simple $St-$coalgebra
$(k[St]^c)^*$ is isomorphic to the $St-$coalgebra
$(k[St]_{\tilde J},\Delta_{\tilde
J})$ for
some unique (up to gauge equivalence) quasitwist
$\tilde J\in
k[St]\ot k[St].$

\begin{Proposition} {\bf [Mov, Propositions
13,14]}\label{movprop1314}
$\tilde J$ is in fact a twist for $k[St]$ (i.e. it is
invertible).
Furthermore, $J$ is the image of $\tilde J$ under the
coalgebra embedding $(k[St]^c)^*\hookrightarrow A_J.$ 
\end{Proposition}
\proof
We only reproduce here the proof of the invertibility of
$\tilde J$ (in a slightly expanded form). Set $C:=k[St]_{\tilde J}.$
Suppose on the contrary that $\tilde J$
is not invertible. Then there exists $0\ne L\in
C^*\ot C^*$ such that $\tilde J L=0.$ Let
$F:C\ot C\raro C\ot C$ be defined by $F(x)=xL.$ Clearly,
$F$ is a morphism of $St\times
St-$representations, and $F\circ \Delta_{\tilde J}=0.$
Thus
the image $\text{Im}(F^*)$ of the morphism of $St\times
St-$representations $F^*:C^*\ot C^*\raro C^*\ot C^*$
is
contained in the kernel of the multiplication map
$m:=(\Delta_{\tilde J})^*.$ Let
$U:=(C^*\ot 1)\text{Im}(F^*)(1\ot C^*).$ Clearly, $U$ is
contained in the kernel of $m$ too.
But, for any $x\in U$ and $g\in St,$ $(1\ot
X_g)x(1\ot X_g)^{-1}\in U.$ Thus, $U$ is a left
$C^*\ot C^*-$module under left multiplication.
Similarly, it is a right module over this algebra under
right multiplication. So, it is a bimodule over
$C^*\ot C^*.$ Since $U\ne 0$, this implies that
$U=C^*\ot C^*.$ This is a contradiction, since we get
that $m=0.$ Hence $\tilde J$ is invertible as desired.
\qed
\begin{Remark} {\rm In the paper [Mov] it is assumed that
the
characteristic of $k$ is equal to $0,$ but all the results
generalize in
a straightforward way to the case when the characteristic
of $k$ is
positive and relatively prime to the order of the
group $G.$}
\end{Remark}

\section{The classification of triangular semisimple
and\\ cosemisimple Hopf algebras}
In this section we describe the classification of
triangular semisimple and
cosemisimple Hopf algebras over {\em any}
algebraically closed field $k,$ given in [EG4].

\subsection{Construction of triangular semisimple and
cosemisimple Hopf algebras from group-theoretical data}
Let $H$ be a finite group
such that $|H|$ is not divisible by the
characteristic of $k.$
Suppose that $V$ is an irreducible projective
representation of $H$ over $k$ satisfying
$\dim(V)=|H|^{1/2}.$ Let $\pi :H\raro \text{PGL}(V)$ be
the projective action of $H$ on $V$, and let
$\tilde \pi:H\to \text{SL}(V)$ be any lifting of this
action
($\tilde\pi$ need not be a homomorphism). We have
$\tilde\pi(x)\tilde\pi(y)=c(x,y)\tilde\pi(xy)$, where
$c$ is a
2-cocycle of $H$ with coefficients in $k^*.$ This
cocycle is
nondegenerate (see Section 3) and hence the
representation
of $H$ on $\text{End}_{k}(V)$ is
isomorphic to the regular representation of $H$ (see
e.g. Proposition \ref{movprop1112}). By
Propositions
\ref{movprop5} and \ref{movprop1314}, this gives rise to
a twist $J(V)$ for $k[H],$ whose equivalence class is
canonically associated to $(H,V)$.

Now, for any group $G\supseteq H,$ whose order is
relatively prime to the
characteristic of $k,$ define a triangular semisimple
Hopf algebra 
\begin{equation}\label{assf}
F(G,H,V):=\left(k[G]^{J(V)},J(V)_{21}^{-1}J(V)\right). 
\end{equation}
We wish to show that it is also cosemisimple.
\begin{lemma}\label{invol}
The Drinfeld element of the triangular semisimple
Hopf algebra $(A,R):=F(G,H,V)$ equals $1.$
\end{lemma}
\proof The Drinfeld element $u$ is a grouplike
element of $A,$
and for any finite-dimensional $A$-module $V$ one has
$\tr|_V(u)=\dim_{\text{Rep}_k(A)}(V)=\dim (V)$ (since
$\text{Rep}_k(A)$ is equivalent to $\text{Rep}_k(G),$
see e.g. [G5]). In particular, we can set $V$
to be the regular representation,
and find that $\tr|_A(u)=\dim(A)\ne 0$ in $k$. But it
is clear that if
$g$ is a nontrivial grouplike element in any
finite-dimensional Hopf
algebra $A,$ then $\tr|_A(g)=0$. Thus, $u=1.$ \qed
\begin{remark}\label{finfdim}
{\rm Lemma \ref{invol} fails for infinite-dimensional 
{\em cotriangular} Hopf algebras, which shows that this
lemma can
not be proved by an explicit computation. We refer the reader
to [EG6] for the study of infinite-dimensional cotriangular
Hopf algebras which are obtained from twisting in function
algebras of affine proalgebraic groups.}
\end{remark}

\begin{corollary} The triangular semisimple Hopf algebra
$(A,R):=F(G,H,V)$ is also cosemisimple.
\end{corollary}
\proof Since $u=1,$ one has $S^2=Id$ and hence $A$ is
cosemisimple
(as $\dim(A)\ne 0$). \qed

Thus we have assigned a triangular semisimple and
cosemisimple Hopf
algebra with Drinfeld element $u=1$ to any triple
$(G,H,V)$ as above.

\subsection{The classification in characteristic $0$}
In this subsection we assume that $k$ is of
characteristic $0.$ We first recall Theorem 2.1 from [EG1] and 
Theorem 3.1 from
[EG4], and state them in a single theorem which is the
key structure theorem for 
triangular semisimple Hopf algebras over $k.$ 

\begin{theorem}\label{mainkey}
Let $(A,R)$ be a triangular semisimple Hopf algebra over an 
algebraically closed field $k$ of characteristic $0,$ with
Drinfeld element $u.$ Set $\tilde R:=RR_u.$ Then there
exist a finite group $G,$ a subgroup $H\subseteq G$ and a
minimal twist $J\in k[H]\ot k[H]$ such that $(A,\tilde
R)$ and $(k[G]^J,J_{21}^{-1}J)$ are isomorphic as triangular
Hopf algebras. Moreover, the data $(G,H,J)$ is unique up to
isomorphism of groups and gauge equivalence of
twists. That is,
if there exist a finite group $G',$ a subgroup
$H'\subseteq G'$ 
and a minimal twist $J'\in k[H']\ot k[H']$ such that
$(A,\tilde
R)$ and $(k[G']^{J'},J_{21}^{'-1}J')$ are isomorphic
as
triangular
Hopf algebras, then there exists an isomorphism of
groups $\phi:G
\raro G'$ such that $\phi(H)=H'$ and $(\phi\ot \phi)(J)$
and $J'$ 
are gauge equivalent as twists for $k[H'].$
\end{theorem}

The proof of this theorem relies, among other things, on 
the following (special case of a) deep
theorem of Deligne on Tannakian categories [De1] in an 
essential way. 

\begin{theorem}\label{cordel}
Let $k$ be an algebraically closed field of characteristic $0,$
and $({\cal C},\ot, {\bf 1},a,l,r,c)$ a $k-$linear abelian
symmetric rigid category with $\text {End}({\bf 1})=k,$
which is semisimple with
finitely many irreducible objects. If categorical
dimensions of objects are
nonnegative integers, then there exist a finite group
$G$ and an equivalence of symmetric rigid categories  
$F:{\cal C}\raro \text {Rep}_k(G).$
\end{theorem}

We refer the reader to [G5, Theorems 5.3, 6.1, Corollary 6.3] for 
a complete and detailed proof of [EG1, Theorem 2.1] and
[EG4, Theorem 3.1], along with a discussion of Tannakian 
categories.

Let $(A,R)$ be a triangular semisimple Hopf algebra
over $k$ whose Drinfeld element $u$ is $1,$ and let
$(G,H,J)$ be the associated group-theoretic
data given in Theorem \ref{mainkey}. 
\begin{proposition}\label{movmin}
The $H-$coalgebra $(k[H]_J,\Delta_J)$ (see
(\ref{mjr})) is simple, and is isomorphic to the
regular representation of $H$ as an $H-$module.
\end{proposition}
\proof By Proposition \ref{movprop1314}, $J$ is the
image of $\tilde J$ under the embedding $k[St]_{\tilde
J}\hookrightarrow k[H]_J.$ Since $J$ is minimal and
$St\subseteq H,$ it follows that $St=H.$ Hence the
result follows from the discussion preceding
Proposition \ref{movprop1314}. \qed

We are now ready to prove our first classification result.
\begin{theorem}\label{class1}
The assignment $F:(G,H,V)\mapsto (A,R)$
is a bijection between:
\ben
\item isomorphism classes of triples
$(G,H,V)$ where $G$ is
a finite group, $H$ is a subgroup of $G$, and
$V$ is an irreducible projective representation of $H$
over $k$ satisfying $\dim(V)=|H|^{1/2}$, and 
\item isomorphism classes of
triangular semisimple Hopf algebras over $k$ with Drinfeld
element $u=1$.
\een
\end{theorem}
\proof We need to construct an assignment $F'$ in the other
direction, and
check that both $F'\circ F$ and $F\circ F'$ are the identity
assignments.

Let $(A,R)$ be a triangular semisimple Hopf algebra over
$k$ whose Drinfeld element $u$ is $1,$ and let $(G,H,J)$
be the associated group-theoretic
data given in Theorem \ref{mainkey}. By Proposition
\ref{movmin}, the $H-$algebra $(k[H]_J)^*$ is simple.
So we see that
$(k[H]_J)^*$ is isomorphic to $\text{End}_{k}(V)$ for
some vector space
$V$, and we have a homomorphism $\pi:H\to \text{PGL}(V)$.
Thus $V$ is a projective representation of $H$.
By Proposition \ref{movprop89}, this representation is
irreducible, and it is obvious that $\dim(V)=|H|^{1/2}.$

It is clear that the isomorphism class
of the representation $V$ does not change if $J$ is
replaced by a twist $J'$ which is gauge
equivalent to $J$ as twists for $k[H].$
Thus, to any isomorphism class of triangular semisimple
Hopf algebras $(A,R)$ over $k,$ with
Drinfeld element $1,$ we have assigned an isomorphism
class of triples $(G,H,V)$. Let us write this as 
\begin{equation}\label{assf'}
F'(A,R):=(G,H,V).
\end{equation}

The identity $F\circ F'=id$ follows from Proposition
\ref{movprop5}.
Indeed, start with $(A,R)\cong
(k[G]^J,J_{21}^{-1}J),$
where $J$ is a minimal twist for $k[H],$ $H$ a subgroup
of $G.$ Then by Proposition \ref{movmin}, 
we have that $(k[H]_J,\Delta_J)$ is a simple
$H-$coalgebra which is 
isomorphic to the regular representation of $H.$ Now let
$V$ be the 
associated irreducible projective representation of $H,$
and $J(V)$ 
the associated twist as in (\ref{assf}). Then 
$(k[H]_J,\Delta_J)$ and
$(k[H]_{J(V)},\Delta_{J(V)})$ are 
isomorphic as $H-$coalgebras, and the claim follows from
Proposition \ref{movprop5}.

The identity $F'\circ F=id$ follows from the
uniqueness part of Theorem \ref{mainkey}. \qed
\begin{remark}\label{hjmin}
{\rm Observe that it follows from Theorem \ref{class1}
that the twist $J(V)$ associated to $(H,V)$ is
{\em minimal}.}
\end{remark}

Now let $(G,H,V,u)$ be a quadruple, in which $(G,H,V)$
is as above,
and $u$ is a central element of $G$ of order $\le 2$.
We extend the map $F$ to quadruples by setting
\begin{equation}\label{assf2}
F(G,H,V,u):=(A,RR_u)\;\text{where}\;(A,R):=F(G,H,V).
\end{equation}
\begin{theorem}\label{clas}
The assignment $F$ given in (\ref{assf2}) is a
bijection between:
\ben
\item isomorphism classes of
quadruples $(G,H,V,u)$ where $G$
is a finite group,
$H$ is a subgroup of $G$, $V$ is an irreducible
projective representation of $H$ over $k$ satisfying
$\dim(V)=|H|^{1/2},$ and $u\in G$ is
a central element of order $\le 2,$
and 
\item isomorphism classes of triangular
semisimple Hopf algebras over $k.$
\een
\end{theorem}
\proof Define $F'$ by $F'(A,R):=(F'(A,RR_u),u),$ where 
$F'(A,RR_u)$ is defined in (\ref{assf'}).
Using Theorem \ref{class1}, it is straightforward to
see that both $F'\circ F$ and
$F\circ F'$ are the identity assignments. \qed

Theorem \ref{clas} implies the following
classification
result for {\em minimal} triangular semisimple Hopf
algebras over $k.$
\begin{proposition}\label{minclas} $F(G,H,V,u)$ is minimal
if and only if $G$ is generated by $H$ and $u$.
\end{proposition}
\proof As we have already pointed out in Remark
\ref{hjmin}, if $(A,R):=F(G,H,V)$
then
the sub Hopf algebra $k[H]^J\subseteq A$ is minimal triangular.
Therefore, if $u=1$ then $F(G,H,V)$ is minimal if and only if
$G=H$. This
obviously
remains true for $F(G,H,V,u)$ if $u\ne 1$ but $u\in H$.
If $u\notin H$
then it is clear that the $R-$matrix of $F(G,H,V,u)$
generates $k[H']$,
where $H'=H\cup uH$. This proves the proposition. \qed
\begin{remark} {\rm As was pointed out already by Movshev, the
theory developed
in [Mov] and extended in [EG4] is an analogue, for finite
groups, of
the theory of quantization of skew-symmetric solutions of
the classical
Yang-Baxter equation, developed by Drinfeld [Dr3]. In
particular,
the operation $F$ is the analogue of the operation of
quantization
in [Dr3].}
\end{remark}

\subsection{The classification in positive characteristic}
In this subsection we assume that $k$ is of
positive characteristic $p,$ and prove an analogue of
Theorem \ref{clas}
by using this theorem itself and the lifting techniques
from [EG5].

We first recall some notation from [EG5].
Let ${\cal O}:=W(k)$ be the ring of Witt vectors 
of $k$ (see e.g. [Se, Sections 2.5, 2.6]), and $K$ the 
field of
fractions of ${\cal O}.$
Recall that ${\cal O}$ is a local complete discrete
valuation ring, and
that the characteristic of $K$ is zero. Let ${\bf m}$
be the maximal ideal in ${\cal O},$ which is
generated by $p.$ One has
${\bf m}^n/{\bf m}^{n+1}=k$ for any $n\ge 0$ (here
${\bf m}^0:={\cal O}$). 

Let $F$ be the assignment defined in (\ref{assf2}).
We now have the
following classification result.
\begin{theorem}\label{clasp}
The assignment $F$ is a bijection between:
\ben
\item isomorphism classes of
quadruples $(G,H,V,u)$ where $G$ is a finite group
of order prime to $p$, $H$ is a subgroup of
$G$, $V$ is an irreducible projective representation of
$H$ over $k$ satisfying
$\dim(V)=|H|^{1/2},$ and $u\in G$ is a central element
of order $\le 2,$ and 
\item isomorphism classes of
triangular semisimple and cosemisimple Hopf
algebras over $k.$
\een
\end{theorem}
\proof As in the proof of Theorem \ref{clas} we need
to construct the
assignment $F'.$ 

Let $(A,R)$ be a triangular semisimple and
cosemisimple Hopf
algebra over $k.$ Lift it (see [EG5]) to a triangular
semisimple
Hopf algebra
$(\bar A,\bar R)$ over $K.$ By Theorem \ref{clas}, we
have that $(\bar 
A\ot_{K}\bar K,\bar R)=F(G,H,V,u).$
We can now reduce $V$ ``$mod\,p$'' to get $V_p$ which
is an irreducible
projective representation of $H$ over the field $k.$
This can be done
since $V$ is defined by a nondegenerate $2-$cocycle
$c$  (see Section 3)
with values in roots of unity of degree $|H|^{1/2}$
(as the only irreducible representation of the simple
$H$-algebra
with basis $\{X_h|h\in H\}$, and relations
$X_gX_h=c(g,h)X_{gh}$).
This cocycle can be reduced $mod\,p$ and remains
nondegenerate
(since the groups of roots of unity of order
$|H|^{1/2}$ in $k$ and $K$
are naturally isomorphic),
so it defines an irreducible projective
representation $V_p.$ Define
$F'(A,R):=(G,H,V_p,u).$ It is shown like in
characteristic $0$ that
 $F\circ F'$ and $F'\circ F$ are the
identity assignments. \qed

The following is the analogue of Theorem \ref{mainkey} 
in positive characteristic.
\begin{corollary}\label{mainkey2}
Let $(A,R)$ be a triangular semisimple and cosemisimple
Hopf algebra over any algebraically closed field $k,$ with
Drinfeld element $u.$ Set $\tilde R:=RR_u.$ Then there
exist a finite group $G,$ a subgroup $H\subseteq G$ and a
minimal twist $J\in k[H]\ot k[H]$ such that $(A,\tilde
R)$ and $(k[G]^J,J_{21}^{-1}J)$ are isomorphic as triangular
Hopf algebras. Moreover, the data $(G,H,J)$ is unique up to
isomorphism of groups and gauge equivalence of
twists. 
\end{corollary}
\begin{proposition}\label{genp}
Proposition \ref{minclas} holds in positive
characteristic as well.
\end{proposition}
\proof As before,
if $(A,R):=F(G,H,V),$ then
the sub Hopf algebra $k[H]^J\subseteq A$ is minimal
triangular.
This follows from the facts that it is true in
characteristic $0,$
and that the rank of a triangular structure
does not change under lifting. Thus, Proposition
\ref{minclas} holds in characteristic $p$. \qed

\subsection{The solvability of the group underlying a minimal
triangular semisimple Hopf algebra}
In this subsection we consider finite groups which admit
a minimal twist as studied in [EG4]. We also consider the 
existence of nontrivial grouplike
elements in triangular semisimple and cosemisimple 
Hopf algebras, following [EG4].

A classical fact about complex representations of finite
groups is that
the dimension of any irreducible representation of a
finite group
$K$ does
not exceed $|K:Z(K)|^{1/2}$, where $Z(K)$ is the center
of $K$.
Groups of central type are those groups for which this
inequality
is in fact an equality. More precisely, a finite group $K$ is
said to be of
{\em central type} if it has an
irreducible representation $V$ such that
$\dim(V)^2=|K:Z(K)|$ (see
e.g. [HI]).
We shall need the following theorem 
(conjectured by Iwahori and Matsumoto in 1964)
whose proof uses the classification of
finite simple groups.
\begin{theorem} {\bf [HI, Theorem 7.3]}\label{ct}
Any group of central type is solvable.
\end{theorem}

As corollaries, we have the following results.
\begin{corollary}\label{solv}
Let $H$ be a finite group which admits a minimal twist. 
Then $H$ is solvable.
\end{corollary}
\proof
We may assume that $k$ has characteristic $0$
(otherwise we can lift to characteristic $0$).
As we showed in the proof of Theorem \ref{class1},
$H$ has an irreducible projective representation $V$ with
$\dim(V)=|H|^{1/2}.$ Let $K$ be a finite central extension
of $H$ with central subgroup $Z,$ such that $V$ lifts to a 
linear representation of $K.$
We have $\text{dim}(V)^2=|K:Z|$.
Since $\text{dim}(V)^2\le |K:Z(K)|$ we get that $Z=Z(K)$
and hence that $K$ is
a group of central type. But by Theorem \ref{ct}, $K$ is
solvable and
hence $H\cong K/Z(K)$ is solvable as well. \qed
\begin{remark}\label{new} {\rm Movshev conjectures in the
introduction to [Mov] that any finite group with a
nondegenerate $2-$cocycle is solvable. As explained in the
Proof of Corollary \ref{solv}, this result follows from
Theorem \ref{ct}.}
\end{remark}
\begin{corollary}\label{ntgl}
Let $A$ be a triangular semisimple and cosemisimple Hopf
algebra over $k$
of dimension bigger than $1$. Then $A$ has a nontrivial
grouplike element.
\end{corollary}
\proof
We can assume that the Drinfeld element $u$ is equal
to $1$ and that $A$ is not cocommutative.
Let $A_m$ be the minimal part of $A$. By Corollary
\ref{solv},
$A_m=k[H]^J$ for a solvable group $H$, $|H|>1$.
Therefore,
$A_m$ has nontrivial $1-$dimensional representations.
Since
$A_m\cong A_m^{*op}$ as Hopf algebras, we get that
$A_m$,
and hence $A$, has nontrivial grouplike elements.\qed

\subsection{Biperfect quasitriangular semisimple Hopf 
algebras}
Corollary \ref{ntgl} motivates the following question.
Let $(A,R)$ be a {\em quasitriangular} semisimple Hopf 
algebra over $k$ with characteristic $0$ (e.g. the
quantum double of a semisimple Hopf algebra), and let 
$\dim(A)>1$. Is it true that $A$ possesses a nontrivial 
grouplike element? We now follow [EGGS] and show that the 
answer to this question is {\em negative}.

Let $G$ be a finite group.
If $G_1$ and $G_2$ are subgroups of $G$ such that
$G=G_1G_2$ and $G_1 \cap G_2 = 1$, we say that
$G=G_1G_2$ is an {\em exact factorization}. In this case
$G_1$ can be identified with $G/G_2,$ and $G_2$ can be
identified with $G/G_1$ as sets, so $G_1$ is a
$G_2-$set and $G_2$ is a $G_1-$set. Note that if $G=G_1G_2$
is an exact factorization, then $G=G_2G_1$ is also an exact
factorization by taking the inverse elements.

Following Kac and Takeuchi [KaG,T] one can construct a
semisimple
Hopf algebra from these data as follows. Consider 
the vector space $H:=\C[G_2]^*\ot \C[G_1].$ Introduce
a product on $H$ by:
\begin{equation}\label{alg}
(\varphi\ot a)(\psi\ot b)=\varphi(a\cdot \psi)\ot ab
\end{equation}
for all $\varphi,\psi\in \C[G_2]^*$ and $a,b\in G_1.$ Here
$\cdot$ denotes the associated action of $G_1$
on the algebra $\C[G_2]^*,$ and $\varphi(a\cdot \psi)$ is
the multiplication of $\varphi$ and $a\cdot \psi$ in the
algebra $\C[G_2]^*.$

Identify the vector spaces 
$$H\ot H=(\C[G_2]\ot \C[G_2])^*\ot
(\C[G_1]\ot \C[G_1])=\text{Hom}_{\C}(\C[G_2]\ot
\C[G_2], \C[G_1]\ot \C[G_1])$$ 
in the usual way, and introduce a coproduct on $H$ by:
\begin{equation}\label{coalg}
(\gD(\varphi\ot a))(b\ot c)=\varphi(bc)a\ot b^{-1}\cdot a
\end{equation}
for all $\varphi\in \C[G_2]^*,$ $a\in G_1$ and $b,c\in G_2.$
Here $\cdot$ denotes the action of $G_2$ on $G_1.$ 

Introduce a counit on $H$ by:
\begin{equation}\label{counit}
\varepsilon(\varphi\ot a)=\varphi(1_G)
\end{equation}
for all $\varphi\in \C[G_2]^*$ and $a\in G_1.$

Finally, identify the vector spaces $H=\C[G_2]^*\ot \C[G_1]=
\text{Hom}_{\C}(\C[G_2],\C[G_1])$ in the usual way, and introduce 
an antipode on $H$ by:
\begin{equation}\label{antipode}
S\left (\sum_{a\in G_1}\varphi_a\ot a \right)(x)=
\sum_{a\in G_1}\varphi_{_{(x^{-1}\cdot a)^{-1}}}\left ((a^{-1}\cdot x)^{-1}
\right)a
\end{equation}
for all $\sum_{a\in G_1}\varphi_a\ot a\in H$ and $x\in G_2,$
where the first $\cdot$ denotes the action of $G_2$ on $G_1,$
and the second one denotes the action of $G_1$ on $G_2.$

\begin{theorem} {\bf [KaG,T]}\label{bicross}
The multiplication, comultiplication, counit and antipode
described in (\ref{alg})-(\ref{antipode}) determine a 
semisimple Hopf algebra structure on
the vector space $H:=\C[G_2]^*\ot \C[G_1].$ 
\end{theorem}

The Hopf algebra $H$ is called the {\em bicrossproduct} Hopf
algebra associated with $G,G_1,G_2,$ and is denoted by
$H(G,G_1,G_2).$
\begin{theorem} {\bf [Ma2]}\label{m}
$H(G,G_2,G_1)\cong H(G,G_1,G_2)^*$ as Hopf algebras.
\end{theorem}

Let us call a Hopf algebra $H$ {\em biperfect} if the groups
$G(H),$ $G(H^*)$ are both trivial. We are ready now to prove:
\begin{theorem}\label{biperf}
$H(G,G_1,G_2)$ is biperfect if and only if $G_1,G_2$ are self
normalizing perfect subgroups of $G.$
\end{theorem}
\proof
It is well known that the category of finite-dimensional
representations of \linebreak $H(G,G_1,G_2)$ is equivalent to
the
category of
$G_1-$equivariant vector bundles on $G_2,$ and hence that the
irreducible representations of $H(G,G_1,G_2)$ are indexed
by pairs $(V,x)$ where $x$ is a representative of a
$G_1-$orbit in $G_2,$ and $V$ is an irreducible
representation
of $(G_1)_x,$ where $(G_1)_x$ is the isotropy subgroup of
$x.$ Moreover, the dimension of the corresponding
irreducible  
representation is
\( \displaystyle \frac{\dim(V)|G_1|}{|(G_1)_x|}.\)
Thus, the
$1-$dimensional representations of $H(G,G_1,G_2)$ are
indexed
by pairs
$(V,x)$ where $x$ is a fixed point of $G_1$ on $G_2=G/G_1$
(i.e. $x\in N_G(G_1)/G_1$), and $V$ is a $1-$dimensional
representation of $G_1.$ The result follows now using
Theorem \ref{m}. \qed

By Theorem \ref{biperf}, in order to construct an
example of a biperfect semisimple Hopf algebra, it
remains to
find a finite group $G$ which
admits an exact factorization $G=G_1G_2,$ where 
$G_1,G_2$ are self normalizing perfect subgroups of $G.$
Amazingly the Mathieu simple group $G:=M_{24}$ of degree $24$
provides such an example!

\begin{theorem}\label{ex}
The group $G$ contains a subgroup
$G_1\cong \text{PSL}(2,23),$ and a subgroup
\linebreak $G_2\cong
A_{7}\ltimes (\Z_2)^4$
where $A_7$ acts on $(\Z_2)^4$ via
the embedding $A_7 \subset
A_8=\text{SL}(4,2)=\text{Aut}((\Z_2)^4)$ (see [AT]).
These subgroups are perfect, self normalizing and $G$
admits
an exact factorization $G=G_1G_2.$ In particular,
$H(G,G_1,G_2)$ is biperfect.
\end{theorem}

We suspect that not only is
$M_{24}$ the smallest example but it may be 
the only finite simple group with a factorization
with all the needed properties.

Clearly, the Drinfeld double $D(H(G,G_1,G_2))$ is an
example of a biperfect quasitriangular semisimple Hopf
algebra.

\subsection{Minimal triangular Hopf algebras
constructed from a bijective 1-cocycle}
In this subsection we describe an explicit way of
constructing
minimal twists for certain solvable groups (hence of
central type groups) given in [EG2]. For
simplicity we let $k:=\C.$
\begin{definition}\label{bij1coc}
Let $G,A$ be finite groups and $\rho :G\to
\text{Aut}(A)$ a homomorphism. 
By a 1-cocycle of $G$ with coefficients in $A$ we
mean a map 
$\pi:G\to A$ which satisfies the equation 
\begin{equation}\label{bij1coc1}
\pi(gg')=\pi(g)(g\cdot \pi(g')),\ g,g'\in G,
\end{equation} 
where $\rho(g)(x)=g\cdot x$ for $g\in G,x\in A.$
\end{definition}

We will be interested in the case when $\pi$ is {\it
a bijection} (so in particular, $|G|=|A|$), 
because of the following proposition. 
\begin{proposition}\label{semitw} Let $G,A$ be finite groups,
$\pi:G\to A$ a 
bijective 1-cocycle, and $J$ a twist for $\C[A]$ which is
$G$-invariant. Then $\bar J:=(\pi^{-1}\ot \pi^{-1})(J)$ 
is a quasitwist for $\C[G].$ 
\end{proposition}
\proof It is obvious that the second equation of
(\ref{t1}) is
satisfied for $\bar J$. So we only have to prove the  
first equation of (\ref{t1}) for $\bar J.$ Write $J=\sum
a_{xy}x\ot y$. Then 
\begin{eqnarray*}\label{stam}
\lefteqn {(\pi\ot\pi\ot\pi)((\Delta\ot Id)(\bar J)(\bar
J\ot 1))}\\
= & & \sum_{x,y,z,t\in
A}a_{xy}a_{zt}\pi(\pi^{-1}(x)\pi^{-1}(z))\ot
\pi(\pi^{-1}(x)\pi^{-1}(t))\ot \pi(\pi^{-1}(y))\\
=& & \sum_{x,y,z,t\in A}a_{xy}a_{zt}x(\pi^{-1}(x)z)\ot
x(\pi^{-1}(x)t)\ot y.
\end{eqnarray*}
Using the $G$-invariance of $J,$ we can remove the
$\pi^{-1}(x)$ 
in the last expression and get 
\begin{equation}\label{stam1}
(\pi\ot\pi\ot\pi)((\Delta\ot Id)(\bar J)(\bar J\ot 1))=
(\Delta\ot Id)(J)(J\ot 1).
\end{equation}
Similarly, 
\begin{equation}\label{stam2}
(\pi\ot\pi\ot\pi)((Id\ot \Delta)(\bar J)(1\ot \bar
J))=
(Id\ot \Delta)(J)(1\ot J).
\end{equation}
But $J$ is a twist, so the right hand sides of
(\ref{stam1}) and
(\ref{stam2}) are equal. 
Since $\pi$ is bijective, this implies equation
(\ref{t1}) for $\bar
J$. \qed
 
Now, given a quadruple $(G,A,\rho,\pi)$ as above
such that $A$ is {\em abelian}, define $\tilde
G:=G\ltimes A^\vee ,$
$\tilde A:=A\times A^\vee ,$ 
$\tilde \rho:\tilde G\raro \text{Aut}(\tilde A)$ by
$\tilde \rho(g)=\rho(g)\times\rho^*(g)^{-1},$ 
and $\tilde \pi:\tilde G\raro \tilde A$ by 
$\tilde \pi(a^*g)=\pi(g)a^*$ for $a^*\in A^\vee $, $g\in G$. 
It is straightforward to 
check that $\tilde \pi$ is a bijective 1-cocycle. We
call the
quadruple 
$(\tilde G,\tilde A,\tilde \rho,\tilde \pi)$ the {\em
double} of
$(G,A,\rho,\pi).$

Consider the element $J\in \C[\tilde A]\ot \C[\tilde A]$ 
given by 
$$
J:=|A|^{-1}\sum_{x\in A,y^*\in A^\vee }e^{(x,y^*)}x\ot y^*,
$$
where $(,)$ is the duality pairing between $A$ and $A^\vee.$
It is straightforward to check that $J$ is a twist
for $\C[\tilde A],$ 
and that it is 
$G$-invariant. This allows to construct the
corresponding element 
\begin{equation}\label{jbar}
\bar J:=
|A|^{-1}\sum_{x\in A,y^*\in A^\vee }
e^{(x,y^*)}\pi^{-1}(x)\ot y^*\in \C[\tilde G]\ot
\C[\tilde G].
\end{equation}

\begin{proposition}\label{1coctw} $\bar J$ is a twist for
$\C[\tilde G],$ and 
\begin{equation}\label{jbarinv}
\bar J^{-1}=|A|^{-1}
\sum_{z\in A,t^*\in A^\vee }e^{-(z,t^*)}\pi^{-1}(T(z))\ot t^*,
\end{equation}
where $T:A\to A$ is a bijective map (not a
homomorphism, in general) 
defined by $$\pi^{-1}(x^{-1})\pi^{-1}(T(x))=1.$$ 
\end{proposition}
\proof Denote the right hand side of (\ref{jbarinv})
by $J'.$ 
We need to check that $J'=\bar J^{-1}$. It is enough
to check it 
after evaluating any $a\in A$ on the second
component of
both sides. We have 
\begin{eqnarray*}
\lefteqn{(1\ot a)({\bar
J})=|A|^{-1}\sum_{x,y^*}e^{(xa,y^*)}\pi^{-1}(x)}\\
= & & \pi^{-1}(a^{-1})(1\ot a)(J')= 
|A|^{-1}\sum_{x,y^*}e^{-(xa^{-1},y^*)}\pi^{-1}(T(x))=
\pi^{-1}(T(a)).
\end{eqnarray*}
This concludes the proof of the proposition. \qed

We can now prove:

\begin{theorem}\label{1cocmintw} Let $\bar J$ be as in
(\ref{jbar}). Then $\bar J$ is
a minimal twist for $\C[\tilde G],$ and
it gives rise to a minimal triangular semisimple Hopf algebra
$(\C[\tilde G]^{\bar J},R^{{\bar J}}),$ with universal
R-matrix 
$$
R^{{\bar J}}=|A|^{-2}\sum_{x,y\in A,x^*,y^*\in
A^\vee }e^{(x,y^*)-(y,x^*)}
x^*\pi^{-1}(x)\ot \pi^{-1}(T(y))y^*.
$$
\end{theorem}
\proof The minimality of ${\bar J}$ follows from the
fact 
that $\{x^*\pi^{-1}(x)|x^*\in A^\vee ,x\in A\}$ and 
$\{\pi^{-1}(T(y))y^*|y\in A,y^*\in A^\vee \}$ are bases of 
$\C[\tilde G],$ and the fact that the matrix
$c_{xx^*,yy^*}=
e^{(x,y^*)-(y,x^*)}$ is invertible (because it is
proportional to
the matrix of Fourier transform on $A\times A^\vee$). \qed

\begin{remark}\label{pavel}
{\rm By Theorem \ref{1cocmintw}, every bijective
1-cocycle $\pi :G\to A$
gives rise to a minimal 
triangular structure on $\C[G\ltimes A^\vee ]$. So it
remains to construct a supply 
of bijective 1-cocycles. This was done in [ESS]. 
The theory of bijective 1-cocycles was developed 
in [ESS], because it was found that they correspond to
set-theoretical
solutions of the quantum Yang-Baxter equation. In
particular, many 
constructions of these 1-cocycles were found. We
refer the reader to 
[ESS] for further details.} 
\end{remark}

We now give two examples of nontrivial minimal triangular
semisimple Hopf algebras. The first one has the least 
possible dimension; namely, dimension $16,$ and the second 
one has dimension $36.$
\begin{example}\label{16}
{\rm Let $G:=\Z_2\times \Z_2$ with generators $x,y,$
and $A:=\Z_4$ with generator $a.$ Define an action of $G$ 
on $A$ by letting $x$ act trivially, and $y$ act
as an automorphism via $y\cdot a=a^{-1}.$ 
Eli Aljadeff pointed out to us that the
group $\tilde G:=(\Z_2\times \Z_2)\ltimes \Z_4$ has a
$2-$cocycle $c$ with coefficients in $\C^*,$ such that
the twisted group algebra $\C[\tilde G]^c$ is simple. This
implies that $\tilde G$ has a minimal twist (see Subsection
4.1). We now use Theorem \ref{1cocmintw} to explicitly
construct our example.

Define a bijective $1-$cocycle $\pi:G\raro A$ as 
follows: $\pi(1)=1,$ $\pi(x)=a^2,$ $\pi(y)=a$ and $\pi(xy)=
a^3.$ Then by Theorem 
\ref{1cocmintw}, ${\bf C}[\tilde G]^{\bar J}$ 
is a non-commutative and non-cocommutative minimal
triangular semisimple Hopf algebra of dimension $16.$

We remark that it follows from the classification of
semisimple
Hopf algebras of dimension $16$ [Kash], that 
the Hopf algebra ${\bf C}[\tilde G]^{\bar J}$ constructed
above,
appeared first in [Kash]. However, our triangular structure 
on this Hopf algebra is new. Indeed, Kashina's triangular
structure on this Hopf algebra is not
minimal, since it arises from a twist of a subgroup of
$\tilde G$
which is isomorphic to $\Z_2\times \Z_2.$}
\end{example}
\begin{example}\label{36}
{\rm Let $G:=S_3$ be the permutation group of three
letters, and
$A:=\Z_2\times 
\Z_3.$ Define an action of $G$ on $A$ by
$s(a,b)=(a,(-1)^{sign(s)}b)$ for $s\in
G,$ $a\in \Z_2$ and $b\in \Z_3.$ Define a bijective
$1-$cocycle
$\pi=(\pi _1,\pi _2):G\raro A$ as follows: $\pi
_1(s)=0$ if
$s$ is even and $\pi _1(s)=1$ if $s$ is odd, and $\pi
_2(id)=0,$ $\pi
_2((123))=1,$ $\pi _2((132))=2,$ $\pi _2((12))=2,$
$\pi _2((13))=0$ and 
$\pi_2((23))=1.$
Then by Theorem \ref{1cocmintw},
${\bf C}[\tilde G]^{\bar
J}$ is a
noncommutative and noncocommutative minimal
triangular semisimple
Hopf algebra of dimension $36.$}
\end{example}

We now wish to determine the group-theoretical data
corresponding, under the bijection of the
classification given in Theorem \ref{clas},
to the minimal triangular semisimple Hopf algebras
constructed in Theorem \ref{1cocmintw}.

Let $H:=G\ltimes A^\vee .$ Following Theorem
\ref{1cocmintw}, we can associate
to this data the element
$$
J:=|A|^{-1}\sum_{g\in G,b\in A^\vee }e^{(\pi(g),b) }b\ot g
$$
(for convenience we use the opposite element to the
one we used before).
We proved
that this element is a minimal twist for $\C[H],$ so
$\C[H]^J$ is a
minimal
triangular semisimple Hopf algebra with Drinfeld
element $u=1.$
Now we wish to find the irreducible projective
representation $V$ of
$H$ which corresponds to $\C[H]^J$ under the
correspondence of 
Theorem \ref{clas}.

Let $V:=\text{Fun}(A,\C)$ be the space of $\C$-valued
functions on $A$.
It has a basis $\{\delta_a|a\in A\}$ of
characteristic functions of
points. Define a projective action $\phi$ of $H$ on $V$
by
\begin{equation}\label{projact}
\phi(b)\delta_a=e^{-(a,b)}\delta_a,\;
\phi(g)\delta_a=\delta_{g\cdot
a+\pi(g)}\;{\rm and}\;\phi(bg)=\phi(b)\phi(g)
\end{equation}
for $g\in G$ and $b\in A^\vee .$ It is straightforward to
verify that this is indeed a projective representation.
\begin{proposition}\label{projrep}
The representation $V$ is irreducible, and corresponds
to $\C[H]^J$ under the bijection of the classification
given in Theorem
\ref{clas}.
\end{proposition}
\proof It is enough to show that the $H$-algebras
$(\C[H]_J)^*$ and $\text{End}_{\C}(V)$
are isomorphic.

Let us compute the multiplication in the algebra
$(\C[H]_J)^*.$ We have
\begin{equation}
\Delta_J(bg)=|A|^{-1}\sum_{g'\in G,b'\in A^\vee }
e^{(\pi(g'),b')}b(g\cdot b')g\otimes bgg'.\label{2}
\end{equation}
Let $\{Y_{bg}\}$ be the dual basis of $(\C[H]_J)^*$ to
the basis $\{bg\}$ of
$\C[H]_J.$
Let $*$ denote the multiplication law dual to the coproduct
$\Delta_J.$
Then, dualizing equation (\ref{2}), we have
\begin{equation}\label{du2}
Y_{b_2g_2}*Y_{b_1g_1}=
e^{(\pi(g_1)-\pi(g_2),b_2-b_1)}Y_{b_1g_2}
\end{equation}
for $g_1,$ $g_2\in G$ and $b_1,$ $b_2\in A^\vee $
(here for convenience we write the operations in $A$
and $A^\vee $ additively).
Define $Z_{bg}:=e^{(\pi(g),b)}Y_{bg}$ for $g\in G$
and $b\in A^\vee .$ In the basis
$\{Z_{bg}\}$ the multiplication law in $(\C[H]_J)^*$
is given by
\begin{equation}\label{multlaw}
Z_{b_2g_2}*Z_{b_1g_1}=e^{(\pi(g_1),b_2)}Z_{b_1g_2}.
\end{equation}

Now let us introduce a left action of $(\C[H]_J)^*$ on
$V$. Set
\begin{equation}\label{lact}
Z_{bg}\delta_a:=e^{(a,b)}\delta_{\pi(g)}.\label{5}
\end{equation}
It is straightforward to check, using (\ref{multlaw}),
that (\ref{lact})
is indeed a left action.
It is also straightforward to compute that this action is
$H-$invariant.
Thus, (\ref{lact}) defines an isomorphism
$(\C[H]_J)^*\to
\text{End}_k(V)$ as
$H-$algebras, which proves the proposition. \qed

\section{The representation theory of cotriangular
semisimple and cosemisimple Hopf algebras}
If $(A,R)$ is a minimal triangular Hopf algebra then
$A$ and $A^{*op}$ are isomorphic as Hopf algebras. But
any nontrivial triangular semisimple and
cosemisimple Hopf algebra $A,$
over any algebraically closed field $k,$ which is {\em
not} minimal, gives rise to a new
Hopf algebra $A^*,$ which is also semisimple and
cosemisimple. These are very
interesting semisimple and cosemisimple Hopf
algebras which arise from finite groups, and they are
abundant by the constructions given in
[EG2,EG4] (see Section 4). Generally, the
dual Hopf algebra of a triangular Hopf algebra is
called {\em cotriangular} in the literature.

In this section we explicitly 
describe the representation theory of cotriangular
semisimple and cosemisimple Hopf algebras
$A^*=(k[G]^J)^*$ studied in [EG3], in terms of
representations of some associated groups. As a
corollary we prove that 
Kaplansky's 6th conjecture [Kap] holds for $A^*;$
that is, that the dimension of any irreducible
representation of $A^*$ divides the dimension of $A.$ 

\subsection{The algebras associated with a twist}
Let $A:=k[H]$ be the group algebra of a finite group
$H$ whose order is relatively prime to the
characteristic of $k.$
Let $J\in A\ot A$ be a {\em minimal} twist, and
$A_1:=(A_J,\Delta_J,\varepsilon)$ be as in
(\ref{mjr}).
Similarly, we define the coalgebra
$A_2:=(_JA,_J\Delta,\varepsilon),$ where $_JA=k[H]$ as
vector spaces, and
\begin{equation}\label{mjl}
_J\Delta:A\raro A\ot A,\;_J\Delta(a)=J^{-1}\Delta(a)
\end{equation}
for all $a\in A.$ Note that since $J$ is a twist,
$_J\Delta$ is indeed coassociative. For $h\in H,$ let
$\delta_h:k[H]\raro k$ be the linear map determined by
$\delta_h(h)=1$ and $\delta_h(h')=0$ for $h\ne h'\in H.$
Clearly the set $\{\delta_h|h\in H\}$ forms a linear basis of 
the dual algebras $(A_1)^*$ and $(A_2)^*.$

\begin{theorem}\label{movshev}
The following hold:
\ben
\item $(A_1)^*$ and $(A_2)^*$ are $H-$algebras via
$$ 
\rho_1(h)\delta_y=\delta_{hy},\;\;\rho_2(h)\delta_y=
\delta_{yh^{-1}}$$
respectively.
\item $(A_1)^*\cong (A_2)^{*op}$ as $H-$algebras
(where $H$ acts on
$(A_2)^{*op}$ as it does on $(A_2)^*$).
\item The algebras $(A_1)^*$ and $(A_2)^*$ are simple,
and are isomorphic as
$H-$modules to the regular representation $R_H$ of $H.$
\een
\end{theorem}
\proof The proof of part 1 is straightforward.

The proof of part 3 follows from Proposition
\ref{movmin} and Part 2.

Let us prove part 2. It is straightforward to
verify that $(S\ot S)(J)=(Q\ot
Q)J_{21}^{-1}\Delta(Q)^{-1}$ where $Q$ is as in
(\ref{tant}) (see e.g. (2.17) in [Ma1,
Section 2.3]).
Hence the map $(A_2)^*\raro (A_1)^{*op},$
$\delta_x\mapsto \delta_{S(x)Q^{-1}}$ determines an
$H-$algebra isomorphism. \qed

Since the algebras $(A_1)^*$, $(A_2)^*$ are simple,
the
actions of $H$ on $(A_1)^*,$ $(A_2)^*$
give rise to projective representations $H\raro
\text{PGL}(|H|^{1/2},k).$
We will denote these projective representations by
$V_1,$
$V_2$ (they can be thought of as the simple modules
over $(A_1)^*,$
$(A_2)^*$, with the induced
projective action of $H$). Note that Part 2 of Theorem
\ref{movshev}
implies that $V_1,$ $V_2$ are
dual to each other, hence that $[V_1]=-[V_2].$

\subsection{The main result}
Let $(A,R)$ be a triangular semisimple and
cosemisimple Hopf algebra
over $k,$ with Drinfeld element $u=1,$ and let
$H,$ $G$ and $J$ be as before.
Consider the dual Hopf algebra $A^*$. It has a basis of
$\delta$-functions $\delta_g$.
The first simple but important fact about the structure
of $A^*$ as an
algebra is the following:
\begin{proposition}\label{propo1}
Let $Z$ be a double coset of $H$ in $G,$ and
$(A^*)_Z:=\oplus_{g\in
Z}k\delta_g\subset A^*.$ Then $(A^*)_Z$ is a subalgebra
of $A^*,$ and
$\displaystyle{A^*=\oplus_Z
(A^*)_Z}$ as algebras.
\end{proposition}
\proof Straightforward. \qed

Thus, to study the representation theory of $A^*,$ it is
sufficient to
describe the representations
of $(A^*)_Z$ for any $Z.$

Let $Z$ be a double coset of $H$ in $G$, and let $g\in Z$.
Let $K_g:=H\cap gHg^{-1},$ and define the embeddings
$\theta_1,\theta_2:K_g\to H$
given by $\theta_1(a)=g^{-1}ag$, $\theta_2(a)=a$. Denote
by $W_i$ the
pullback
of the projective $H$-representation $V_i$ to $K_g$ by
means of
$\theta_i$, $i=1,2$.

Our main result is the following theorem, which is
proved in the next
subsection.
\begin{theorem}\label{maincot}
Let $W_1,W_2$ be as above, and let $(\hat
K_g,\tilde\pi_{_W})$ be any
linearization of the projective representation
$W:=W_1\otimes W_2$ of $K_g.$ Let $\zeta$ be the kernel
of the projection
$\hat K_g\to K_g,$ and $\chi:\zeta\to k^*$ be the
character by which
$\zeta$ acts in $W$. Then there exists a 1-1
correspondence between:
\ben
\item isomorphism classes of
irreducible representations of $(A^*)_Z$ and
\item isomorphism classes of
irreducible representations of
$\hat K_g$ with $\zeta$ acting by $\chi$.
\een
Moreover, if a representation $Y$ of $(A^*)_Z$ corresponds 
to a representation
$X$ of $\hat K_g,$ then
$$\dim(Y)=\frac{|H|}{|K_g|}\dim(X).$$ 
\end{theorem}

As a corollary we get Kaplansky's 6th conjecture [Kap] for
cotriangular semisimple and cosemisimple Hopf algebras.
\begin{corollary}\label{kap}
The dimension of any irreducible representation of a
cotriangular semisimple and cosemisimple Hopf algebra
over $k$ divides the dimension of the Hopf algebra.
\end{corollary}
\proof Since $\dim(X)$ divides $|K_g|$ (see e.g. [CR, 
Proposition 11.44]), we have that \linebreak
$\displaystyle{\frac{|G|}{\frac{|H|}{|K_g|}\dim(X)}=
\frac{|G|}{|H|}\frac{|K_g|}{\dim(X)}}$ and the result
follows. \qed

In some cases the classification of representations of
$(A^*)_Z$ is even simpler.
Namely, let $\ol{g}\in \text{Aut}(K_g)$ be given by
$\ol{g}(a)=g^{-1}ag.$ Then we have:
\begin{corollary}\label{corf}
If the cohomology class
$[W_1]$ is $\ol{g}-$invariant then irreducible
representations of
$(A^*)_Z$ correspond in a 1-1
manner to irreducible representations of
$K_g$, and if $Y$ corresponds to $X$ then
$\displaystyle{\dim(Y)=\frac{|H|}{|K_g|}\dim(X)}.$
\end{corollary}
\proof
For any $\ga \in \text{Aut}(K_g)$ and $f\in
\text{Hom}((K_g)^n,k^*),$ let $\ga \circ
f\in \text{Hom}((K_g)^n,k^*)$
be given by $(\ga \circ f)(h_1,\dots,h_n)=f(\ga
(h_1),\dots,\ga (h_n))$
(which determines the
action of $\ga$ on $H^i(K_g,k^*)).$ Then it follows from
the identity
$[V_1]=-[V_2]$
that $[W_1]=-\ol{g}\circ[W_2].$ Thus, in our situation
$[W]=0$,
hence $W$ comes from a linear representation of
$K_g$. Thus, we can set $\hat K_g=K_g$ in the theorem,
and the result
follows.\qed
\begin{example} {\rm Let $k:=\C.$ Let $p>2$ be a prime
number, and
$H:=(\Z/p\Z)^2$
with the standard symplectic
form $(,):H\times H\to k^*$ given by
$\displaystyle{\left((x,y),(x',y')\right)=
e^{2\pi i(xy'-yx')/p}.}$ Then the element
$\displaystyle{J:=p^{-2}\sum_{a,b\in H}(a,b)a\ot b}$ is a
minimal
twist
for $\C[H].$ Let $g\in \text{GL}_2(\Z/p\Z)$ be an
automorphism of
$H$,
and $G_0$ be the cyclic group generated by $g$.
Construct the group $G:=G_0\ltimes H.$ 
It is easy to
see that in this case, the
double cosets are ordinary
cosets $g^kH$, and $K_{g^k}=H$. Moreover, one can
show either explicitly
or using Proposition \ref{movprop89}, that $[W_1]$ is
a generator
of $H^2(H,\C^*)$ which is isomorphic to $\Z/p\Z.$ The
element
$g^k$ acts on $[W_1]$ by multiplication by
$\text{det}(g^k)$.
Therefore,
by Corollary \ref{corf},
the algebra $(A^*)_{g^kH}$ has $p^2$ 1-dimensional
representations
(corresponding to linear representations of $H$) if
$\text{det}(g^k)=1$.

However, if $\text{det}(g^k)\ne 1$, then $[W]$ generates
$H^2(H,\C^*)$. Thus, $W$ comes from a linear
representation of the Heisenberg
group $\hat H$ (a central extension of $H$ by
$\Z/p\Z$) with some
central character $\chi$. Thus, $(A^*)_{g^kH}$ has
one
p-dimensional
irreducible representation, corresponding
to the unique irreducible representation of $\hat H$
with central character $\chi$ (which is $W$).}
\end{example}

\subsection{Proof of Theorem \ref{maincot}}
Let $Z\subset G$ be a double coset of $H$ in $G.$
For any $g\in Z$ define the linear map
$$F_g:(A^*)_Z\raro (A_2)^*\ot
(A_1)^*,\;\;\delta_y\mapsto
\sum_{h,h'\in
H:y=hgh'}\delta_h\ot
\delta_{h'}.$$
\begin{proposition}\label{prop4}
Let $\rho_1,\rho_2$ be as in Theorem \ref{movshev}. Then:
\ben
\item The map $F_g$ is an injective homomorphism of
algebras.
\item $F_{aga'}(\varphi)=(\rho_2(a)\ot
\rho_1(a')^{-1})F_g(\varphi)$ for any
$a,a'\in H,$ $\varphi\in
(A^*)_Z.$
\een
\end{proposition}
\proof
1. It is straightforward to verify that the map
$(F_g)^*:A_2\ot
A_1\raro A_Z$ is determined by $h \ot h'\mapsto hgh',$
and that it is a
surjective homomorphism of
coalgebras. Hence the result follows.

2. Straightforward. \qed

 For any $a\in K_g$ define $\rho(a)\in
\text{Aut}((A_2)^*\ot (A_1)^*)$ by
$\rho(a)=\rho_2(a)\ot \rho_1(a^g),$ where
$a^g:=g^{-1}ag$ and
$\rho_1,\rho_2$ are as in Theorem
\ref{movshev}. Then $\rho$ is an action of $K_g$ on
$(A_2)^*\ot (A_1)^*.$
\begin{proposition}\label{prop5}
Let $U_g:=((A_2)^*\otimes (A_1)^*)^{\rho(K_g)}$ be the
algebra of invariants.
Then $\text{Im}(F_g)=U_g,$
so $(A^*)_Z\cong U_g$ as algebras.
\end{proposition}
\proof
It follows from Proposition \ref{prop4} that
$\text{Im}(F_g)\subseteq U_g,$
and 
$\displaystyle{\text{rk}(F_g)=\dim((A^*)_Z)=
\frac{|H|^2}{|K_g|}}.$ On the
other hand, by Theorem \ref{movshev},
$(A_1)^*,(A_2)^*$ are isomorphic to the regular
representation $R_H$ of $H.$
Thus, $(A_1)^*,(A_2)^*$ are
isomorphic to $\displaystyle{\frac{|H|}{|K_g|}R_{K_g}}$ as
representations of $K_g$, via
$\rho_1(a),\rho_2(a^g)$. Thus,
$\displaystyle{(A_2)^*\ot (A_1)^*\cong
\frac{|H|^2}{|K_g|^2}(R_{K_g}\ot
R_{K_g})\cong \frac{|H|^2}{|K_g|}R_{K_g}}.$ So $U_g$ has
dimension $|H|^2/|K_g|,$ and the result follows. \qed

Now we are in a position
to prove Theorem \ref{maincot}. Since $W_1\ot
W_1^*\cong (A_1)^*$ and $W_2\ot W_2^*\cong (A_2)^*,$
it
follows from Theorem \ref{movshev} that
$\displaystyle{W_1\ot W_2\otimes W_1^*\otimes W_2^*\cong
\frac{|H|^2}{|K_g|}R_{K_g}}$ as
$\hat{K_g}$ modules. Thus, if
$\chi_{_W}$ is the character of $W:=W_1\otimes W_2$ as
a $\hat{K_g}$
module then
$$|\chi_{_W}(x)|^2=0,\,x\notin \zeta\;\;
and\;\;|\chi_{_W}(x)|^2=|H|^2,\,x\in\zeta.$$
Therefore,
$$\chi_{_W}(x)=0,\,x\notin \zeta\;\;
and\;\;\chi_{_W}(x)=|H|\cdot
x_{_W},\,x\in\zeta,$$
where $x_{_W}$ is the root of unity by which $x$ acts
in $W$.
Now, it is clear from the definition of $U_g$ (see
Proposition
\ref{prop5}) that
$U_g=\text{End}_{\hat{K_g}}(W).$
Thus if $\displaystyle{W=\bigoplus_{M\in Irr(\hat{K_g})}
W(M)\ot M},$
where
$W(M):=\text{Hom}_{\hat{K_g}}(M,W)$ is the
multiplicity space,
then
$\displaystyle{U_g=\bigoplus_{M:W(M)\ne
0}\text{End}_k(W(M))}.$
So $\{W(M)|W(M)\ne 0\}$ is the set of
irreducible representations of $U_g.$ Thus the following
implies the
theorem:

\noindent
{\bf Lemma} {\em The following hold:
\ben
\item $W(M)\ne 0$ if and only if for all $x\in \zeta,$
$x_{_{|M}}=x_{_{|W}}.$
\item If $W(M)\ne 0$ then
$\displaystyle{\dim(W(M))=\frac{|H|}{|K_g|}\dim(M)}.$
\een
}

\noindent
{\bf Proof:} The "only if" part of 1 is
clear. For the "if"
part compute $\dim(W(M))$ as
the
inner product $(\chi_{_{W}},\chi_{_M}).$ We have
$$(\chi_{_{W}},\chi_{_M})=\sum_{x\in
\zeta}\frac{|H|}{|\hat{K_g}|}x_{_{|W}}\cdot
\dim(M)\cdot \bar x_{_{|M}}.$$ If
$x_{_{|M}}=x_{_{|W}}$ then
$$(\chi_{_{W}},\chi_{_M})=\sum_{x\in
\zeta}\frac{|H|}{|\hat{K_g}|}\dim(M)=
\frac{|H||\zeta|}{|\hat{K_g}|}\dim(M)
=\frac{|H|}{|K_g|}\dim(M).$$
This proves part 2 as well. \qed

This concludes the proof of the theorem. \qed

\section{The pointed case}
In this section we consider finite-dimensional triangular
pointed Hopf algebras over an algebraically
closed field $k$ of characteristic $0,$ and in particular describe 
the classification and explicit construction of minimal
triangular pointed Hopf algebras, given in [G4].
Throughout the section, unless otherwise specified,
the ground field $k$ will be assumed to be algebraically closed with
characteristic $0.$

\subsection{The antipode of triangular pointed Hopf
algebras} 
In this subsection we prove that the fourth power
of the antipode of any triangular pointed Hopf algebra
$(A,R)$ is the identity. Along the way we prove that
the group algebra of the group of grouplike elements of
$A_R$ (which must be abelian) admits a
minimal triangular structure and consequently that $A$
has the structure of a biproduct [R1].
 
\begin{theorem}\label{gen}
Let $(A,R)$ be a minimal triangular pointed Hopf
algebra over $k$ with Drinfeld 
element $u,$ 
and set $K:=k[\G(A)].$ Then there exists a projection
of Hopf
algebras $\pi:A\raro K,$ and consequently
$A=B\times K$ is a biproduct where $B:=\{x\in A|(I\ot 
\pi)\Delta(x)=x\ot 1\}\subseteq A.$ Moreover, $K$
admits a minimal 
triangular structure with Drinfeld element $u_{_K}=u.$
\end{theorem}
\proof Since $\G(A)$ is abelian, $K^*\cong K$
and $K\cong k[\G(A^{*cop})]$ as Hopf algebras.
Hence,  
$\dim(K^*)=\dim(k[\G(A^{*})])$. Consider the
series of Hopf algebra homomorphisms
$$
K\stackrel{i}{\hookrightarrow}A^{cop}\stackrel{(f_R)^{-1}}
{\longrightarrow}A^*\stackrel{i^*}{\longrightarrow} K^*,
$$ 
where $i$ is the inclusion map.
Since $A^*$ is pointed it follows from the
above remarks that
$i^*_{|k[\G(A^{*})]}:k[\G(A^{*})]\raro 
K^*$ is an isomorphism of Hopf algebras (see e.g. [Mon,
5.3.5] ), and 
hence that $i^*\circ (f_R)^{-1}\circ i$ determines a minimal
quasitriangular 
structure on $K^*.$ This structure is in fact
triangular since
$(f_R)^{-1}$ determines a triangular structure on $A^*.$
Clearly, $(i^*\circ (f_R)^{-1}\circ 
i)(u)=(u_{_{K^*}})^{-1}=u_{_{K^*}}$ is the Drinfeld element
of $K^*.$
Since $K$ and $K^*$ are isomorphic as Hopf
algebras we 
conclude that $K$ admits a minimal triangular
structure with Drinfeld
element $u_{_K}=u.$

Finally, set $\varphi:=i^*\circ (f_R)^{-1}\circ i$ and
$\pi:=\varphi^{-1}\circ i^* 
\circ (f_R)^{-1}.$ Then $\pi:A\raro K$ is onto, and
moreover 
$\pi\circ i=\varphi^{-1}\circ i^*\circ (f_R)^{-1}\circ
i=\varphi^{-1}\circ 
\varphi=id_{K}.$ Hence $\pi$ is a projection of Hopf
algebras and 
by [R1], $A=B\times K$ is a biproduct where $B:=\{x\in
A|(I\ot 
\pi)\Delta(x)=x\ot 1\}$ as desired. This concludes
the proof of the 
theorem. \qed

\begin{theorem}\label{ant}
Let $(A,R)$ be any triangular pointed Hopf algebra
with antipode $S$
and Drinfeld element $u$ over
any field $k$ of characteristic $0.$  
Then $S^4=Id.$ If in
addition $A_m$ is 
not semisimple and $A$ is finite-dimensional then
$\dim(A)$ is divisible by $4.$
\end{theorem}
\proof We may assume that $k$ is algebraically closed.
By (\ref{ssqu}), $S^2(a)=uau^{-1}$ for all 
$a\in A.$ Let $K:=k[\G(A_m)].$ Since $u\in A_m,$ and
by Theorem 
\ref{gen}, $u=u_{_K}$ and $(u_{_K})^2=1,$ we have that
$S^4=Id.$ 

In order to prove the second claim, we may assume
that $(A,R)$ 
is minimal (since by [NZ], $\dim(A_m)$ divides
$\dim(A)$).
Since $A$ is not semisimple it follows from [LR1] that
$S^2\ne Id,$ 
and hence that $u\ne 1.$ In particular, $|\G(A)|$ is
even. Now,
let $B$ be as
in Theorem \ref{gen}. Since $S^2(B)=B,$ $B$ has a basis
$\left \{a_i,b_j|S^2(a_i)=a_i,S^2(b_j)=-b_j,1\le i\le
n,1\le j\le m\right \}.$
Hence by Theorem \ref{gen}, 
$$\left \{a_ig,b_jg|g\in \G(A),1\le i\le n,1\le j\le
m\right \}$$ 
is a basis of $A.$ Since by [R3], $\tr(S^2)=0,$ we
have that 
$0=\tr(S^2)=|\G(A)|(n-m),$ which implies that $n=m,$
and hence that
$\dim(B)$ is even as well. \qed

In fact, the first part of Theorem \ref{ant} can be
generalized.

\begin{theorem}\label{p}
Let $(A,R)$ be a finite-dimensional quasitriangular
Hopf algebra with
antipode $S$ over any field $k$ of characteristic $0,$
and suppose that
the Drinfeld element $u$ of $A$ acts as a scalar in
any irreducible
representation of $A$ (e.g. when $A^*$ is pointed).
Then $u=S(u)$ and in
particular $S^4=Id.$
\end{theorem}
\proof We may assume that $k$ is algebraically
closed. In any
irreducible
representation $V$ of $A,$
$\tr_{|V}(u)=\tr_{|V}(S(u))$ (see
Subsection 2.1). Since $S(u)$ also acts as a scalar in
$V$ (the dual of
$S(u)_{|V}$ equals $u_{|V^*}$) it follows that
$u=S(u)$ in any irreducible representation of $A.$
Therefore,
there exists
a basis of $A$ in which the operators of left
multiplication by 
$u$ and $S(u)$ are represented by upper triangular
matrices with
the same
main diagonal. Hence the special grouplike element
$uS(u)^{-1}$ is
unipotent. Since it has a finite order we conclude that
$uS(u)^{-1}=1,$
and hence that $S^4=Id.$ \qed

\begin{remark}\label{ge}
{\rm If $(A,R)$ is a minimal triangular pointed Hopf
algebra then all its
irreducible representations are $1-$dimensional.
Hence Theorem
\ref{p} is
applicable, and the first part of Theorem \ref{ant}
follows.}
\end{remark}
\begin{example}\label{sweed}
{\rm Let $A$ be Sweedler's $4-$dimensional Hopf
algebra [Sw]. It is generated as
an algebra by a
grouplike element $g$ and a $1:g$ skew primitive
element $x$ satisfying
the relations
$g^2=1,$ $x^2=0$ and $gx=-xg.$ It is known that
$A$ admits minimal
triangular structures all of which with $g$ as the
Drinfeld element [R2]. In
this example, $K=k[<g>]$ and $B=sp\{1,x\}.$ Note that
$g$ is
central in $K$
but is not central in $A,$ so $(S_{|K})^2=Id$ but
$S^2\ne Id$ in
$A.$
However, $S^4=Id.$}
\end{example}

\subsection{Construction of minimal triangular pointed Hopf
algebras}
In this section we give a method for the construction of minimal
triangular pointed Hopf algebras which are {\em not}
necessarily semisimple.

Let $G$ be a finite abelian group, and $F:G\times G\raro k^*$ be a 
non-degenerate skew symmetric bilinear form on $G.$ That is,
$F(xy,z)=F(x,z)F(y,z),$ $F(x,yz)=F(x,y)F(x,z),$ $F(1,x)=F(x,1)=1,$ 
$F(x,y)=F(y,x)^{-1}$ for all $x,y,z\in G,$ and the map
$f:G\raro G^\vee$ defined by $<f(x),y>=F(x,y)$ for
all $x,y\in G$ is an isomorphism. Let
$U_F:G\raro \{-1,1\}$ be defined by $U_F(g)=F(g,g).$ Then $U_F$ is
a homomorphism of groups. Denote $U_F^{-1}(-1)$ by $I_F.$

\begin{definition}\label{datum}
Let $k$ be an algebraically closed field of characteristic zero.
A datum \linebreak
${\cal D}=(G,F,n)$ is a triple where
$G$ is a finite abelian group, $F:G\times G\raro k^*$ 
is a non-degenerate skew symmetric bilinear form on $G,$ and $n$ is a
non-negative integer function $I_F\raro \Z^+,$ $g\mapsto n_g.$
\end{definition}

\begin{remark}\label{form}
1) The map $f:k[G]\raro k[G^\vee]$ determined by $<f(g),h>=F(g,h)$ for 
all $g,$ $h\in G$ determines a minimal triangular structure on $k[G^\vee].$

\noin
2) If $I_F$ is not empty then $G$ has an even order.
\end{remark}

To each datum ${\cal D}$ we associate a Hopf algebra $H({\cal D})$ in the
following way. For each $g\in I_F,$ 
let $V_g$ be a vector space of dimension $n_g,$ and let
${\cal B}=\bigoplus_{g\in I_F} V_g.$
Then $H({\cal D})$ is generated as an algebra by $G\cup {\cal B}$ with the
following additional relations (to those of the group $G$ and the vector
spaces $V_g$'s):
\begin{equation}\label{3p}
xy=F(h,g)yx\;\;and\;\;xa=F(a,g)ax 
\end{equation}
for all $g,h\in I_F,$ $x\in V_g,$ $y\in V_h$ and $a\in G.$

The coalgebra structure of $H({\cal D})$ is determined by letting $a\in G$
be a grouplike element and $x\in V_g$ be a $1:g$ skew primitive element
for all $g\in I_F.$ In particular, $\varepsilon(a)=1$ and
$\varepsilon(x)=0$ for all $a\in G$ and $x\in V_g.$

In the special case where $G=\Z_2=\{1,g\},$ $F(g,g)=-1$ and $n:=n_g,$
the associated Hopf algebra will be denoted by $H(n).$ Clearly,
$H(0)=k\Z_2,$ $H(1)$ is Sweedler's $4-$dimensional Hopf algebra, and
$H(2)$ is the $8-$dimensional Hopf algebra studied in [G1, Section 2.2] in
connection with KRH invariants of knots and $3-$manifolds. We remark that the
Hopf algebras $H(n)$ are studied in [PO1,PO2] where they are denoted by
$E(n).$

For a finite-dimensional vector space $V$ we let $\bigwedge V$ denote the
exterior algebra of $V.$ Set $B:=\bigotimes_{g\in I_F}\bigwedge V_g.$

\begin{proposition}\label{dual}
The following hold:
\ben
\item The Hopf algebra $H({\cal D})$ is pointed and
$\G(H({\cal D}))=G.$

\item $H({\cal D})=B\times k[G]$ is a biproduct. 

\item $H({\cal D})_1= k[G]\bigoplus (k[G]{\cal B}),$ and
$P_{a,b}(H({\cal
D}))=sp\{a-b\}\bigoplus aV_{a^{-1}b}$ for all $a,b\in G$
(here we agree that $V_{a^{-1}b}=0$ if $a^{-1}b\notin I_F$).
\een
\end{proposition}

\proof 
Part 1) follows since (by definition) $H({\cal D})$ is generated as an
algebra by grouplike elements and skew primitive elements.
Now, it is straightforward to verify that the map
$\pi:H({\cal D})\raro k[G]$ determined by
$\pi(a)=a$ and $\pi(x)=0$ for all $a\in G$ and $x\in
{\cal B},$ is a projection of Hopf algebras. Since $B=\{x\in H({\cal
D})|(I\ot
\pi)\Delta(x)=x\ot 1\},$ Part 2) follows from [R1]. Finally,
by Part 2), $B$ is a braided graded Hopf algebra in the Yetter-Drinfeld
category $_{k[G]}^{k[G]}{\cal YD}$ (see e.g. [AS]) with respect to the
grading where the elements of ${\cal B}$ are homogeneous of
degree $1.$ Write $B=\bigoplus_{n\ge 0}B(n),$ where $B(n)$
denotes the homogeneous component of degree $n.$ Then, $B(0)=k1=B_1$
(since $B\cong H({\cal D})/H({\cal D})k[G]^+$ as coalgebras, it is
connected). Furthermore, by similar arguments used in the proof of
[AS, Lemma 3.4], $P(B)=B(1)={\cal B}.$ But
then by [AS, Lemma 2.5], $H({\cal D})$ is coradically graded (where the
$nth$ component $H({\cal D})(n)$ is just $B(n)\times k[G]$) which means
by definition that $H({\cal D})_1= H({\cal D})(0)\bigoplus H({\cal
D})(1)=k[G]\bigoplus (k[G]{\cal B})$ as desired. The second statement
of Part 3) follows now, using (1), by counting dimensions. \qed

In the following we determine {\em all} the minimal triangular structures
on $H({\cal D}).$ Let \linebreak $f:k[G]\raro k[G^\vee]$ be the isomorphism
from
Remark \ref{form} 1), and set $I_F':=\{g\in I_F|n_g\ne 0\}.$ Let
$\Phi$ be the set of all isomorphisms $\varphi:G^\vee\raro G$
satisfying $\varphi^*(\ga)=\varphi(\ga ^{-1})$ for all $\ga\in G^\vee,$ and
$(\varphi\circ f)(g)=g$ for all $g\in I_F'$ (here we identify $G$ with $G^\vee
{^\vee }$).

Extend any $\ga \in G^\vee$ to an algebra homomorphism 
$H({\cal D})\raro k$ by setting $\ga (z)=0$ for all $z\in {\cal B}.$
Extend any $x\in V_g^*$ to $P_x\in H({\cal D})^*$ by setting $<P_x,ay>=0$
for all $a\in G$ and $y\in \bigotimes_{g\in I_F}\bigwedge V_g$ of degree
different from $1,$ and $<P_x,ay>=\delta_{g,h}<x,y>$ for all $a\in G$ and
$y\in V_h.$ We shall identify the vector
spaces $V_g^*$ and $\{P_x|x\in V_g^*\}$ via the map $x\mapsto P_x.$

For $g\in I_F',$ let $S_g(k)$ be the set of all isomorphisms $M_g:V_g^*\raro
V_{g^{-1}}.$ Let $S(k)\subseteq \times _{g\in I_F'} S_g(k)$ be the set of all 
tuples $(M_g)$ satisfying $M_g^*=M_{g^{-1}}$ for all $g\in I_F'.$

\begin{theorem}\label{quas}
1) For each $T:=(\varphi,(M_g))\in \Phi\times S(k),$ there exists a unique Hopf
algebra isomorphism $f_T:H({\cal
D})^{*cop}\raro H({\cal D})$ determined by $\ga \mapsto
\varphi(\ga)$ and $P_x\mapsto M_{g}(x)$ for $\ga\in G^\vee$ and $x\in
V_g^*.$  

\noin
2) There is a one to one correspondence between
$\Phi\times S(k)$ and the
set of minimal triangular structures on $H({\cal D})$ given by $T\mapsto
f_T.$
\end{theorem}

\proof We first show that $f_T$ is a well defined isomorphism of Hopf
algebras.
Using Proposition \ref{dual} 2), it is straightforward to
verify that $\Delta(P_x)=\varepsilon\ot 
P_x+P_x\ot f(g^{-1}),$ $P_x\ga=<\ga,g>\ga P_x,$ and
$P_xP_y=F(h,g)P_yP_x$ for all $\ga\in G^\vee,$  
$g,h\in I_F$, $x\in V_g^*$ and $y\in V_h^*.$
Let ${\cal B}^*:=\{P_x|x\in V_g^*,\;g\in I_F\},$ and $H$ be the
sub Hopf algebra of $H({\cal D})^{*cop}$ generated as an algebra by
$G^\vee\bigcup {\cal B}^*.$
Then, using (4) and our assumptions on $T,$
it is straightforward to verify that the map $f_T^{-1}:H({\cal
D})\raro H$ determined by $a\mapsto \varphi^{-1}(a)$ and
$z\mapsto M_{g}^{-1}(z)$ for $a\in G$ and $z\in
V_{g^{-1}},$ is a surjective homomorphism of Hopf algebras.
Let us verify for instance
that $f_T^{-1}(za)=F(a,g)f_T^{-1}(az).$ Indeed,
this is equivalent to $<\varphi^{-1}(a),g>=<f(a),g>$ which in turn
holds by our assumptions on $\varphi.$ Now, using Proposition 4.3 3), it
is
straightforward to verify that $f_T^{-1}$ is injective on
$P_{a,b}(H({\cal D}))$ for all $a,b\in G.$
Since $H({\cal D})$ is pointed, $f_T^{-1}$ is also injective
(see e.g. [Mon, Corollary 5.4.7]). This implies that
$H=H({\cal
D})^{*cop},$ and that $f_T:H({\cal D})^{*cop}\raro H({\cal D})$ is an
isomorphism of Hopf algebras as desired. Note that in particular,
$G^\vee=G(H({\cal D})^*).$

The fact that $f_T$ satisfies (\ref{qt5}) follows from 
a straightforward computation (using (\ref{3p})) since it is enough to
verify it for algebra generators $p\in G^\vee\cup {\cal B}^*$ of $H({\cal
D})^{*cop},$ and $a\in G\cup {\cal B}$ of $H({\cal D}).$ 

We have to show that $f_T$ satisfies (\ref{t}).
Indeed, it is straightforward to verify that \linebreak $f_T^*:H({\cal
D})^{*op}\raro H({\cal D})$ is determined by $\ga \mapsto
\varphi(\ga^{-1})$ and $P_x\mapsto gM_{g}(x)$ for
$\ga\in 
G^\vee$ and $x\in V_g^*.$ Hence,
$f_T^*=f_T\circ S,$
where $S$ is the antipode of $H({\cal D})^*,$ as desired. 

We now have to show that {\em any} minimal triangular structure on $H({\cal
D})$ comes from $f_T$ for some $T.$ 
Indeed, let ${\bf f}:H({\cal D})^{*cop}\raro H({\cal D})$ 
be any Hopf isomorphism. Then ${\bf f}$ must map
$G^\vee$ onto $G,$ $\{f(g^{-1})|g\in I_F'\}$ onto $I_F',$ and 
$P_{f(g^{-1}),\varepsilon}(H({\cal D})^{*cop})$ bijectively onto
$P_{1,\varphi(f(g^{-1}))}(H({\cal D})).$ Therefore
there exists an invertible operator $M_g:V_{g}^*\raro
V_{\varphi(f(g^{-1}))}$ such that ${\bf f}$ is determined by $\ga\mapsto
\varphi(\ga)$ and $P_x\mapsto M_g(x).$ Suppose
${\bf f}$ satisfies (\ref{qt5}). Then letting $p=P_x$ and $a\in G$ in
(\ref{qt5}) yields that $a{\bf f}(P_x)=F(a,g){\bf f}(P_x)a$ for all $a\in
G.$ But by (\ref{3p}), this is equivalent to $(\varphi\circ f)(g)=g$
for all
$g\in I_F'.$ Since by Theorem \ref{gen}, $\varphi:k[G^\vee]\raro k[G]$
determines a minimal triangular structure on $k[G]$ it follows that
$\varphi\in \Phi.$ Since ${\bf f}:H({\cal D})^{*cop}\raro H({\cal D})$ 
satisfies (\ref{t}), $(M_g)\in S(k),$
and hence ${\bf f}$ is of the form $f_T$ for some $T$ as desired. \qed 

For a triangular structure on $H({\cal D})$ corresponding to the map 
$f_T,$ we let $R_T$ denote the corresponding $R-$matrix. 

\begin{remark}
{\rm Note that if $n_{g^{-1}}\ne n_g$ for some $g\in I_F',$ 
then $S(k)$ is empty and $H({\cal D})$ does {\em not} have
a minimal triangular structure.
}
\end{remark}

\subsection{The classification of minimal triangular pointed
Hopf algebras}
In this subsection we use Theorems \ref{gen}, \ref{ant} and
[AEG, Theorem 6.1] to classify
minimal triangular pointed Hopf algebras. Namely, we prove:

\begin{theorem}\label{mainpoint}
Let $(A,R)$ be a minimal triangular pointed Hopf
algebra over an algebraically closed field $k$ of characteristic
$0.$ There exist a datum ${\cal D}$ and $T\in \Phi\times
S(k)$ such that
$(A,R)\cong (H({\cal D}),R_T)$ as triangular Hopf algebras.
\end{theorem}

Before we prove Theorem \ref{mainpoint} we need to fix some notation
and prove
a few lemmas.

In what follows, $(A,R)$ will always be a minimal triangular
pointed Hopf
algebra over $k,$ $G:=\G(A)$ 
and $K:=k[\G(A)].$ For any $g\in G,$ $P_{1,g}(A)$ is a
$<g>-$module 
under conjugation, and $sp\{1-g\}$ is a submodule of
$P_{1,g}(A).$ Let $V_g\subset P_{1,g}(A)$ be its complement, and set
$n_g:=\dim(V_g).$

By Theorem \ref{gen}, $A=B\times K$ where $B=\{x\in A|(I\ot 
\pi)\Delta(x)=x\ot 1\}\subseteq A$ is a 
left coideal subalgebra of $A$ (equivalently, $B$ is an object in the
Yetter-Drinfeld category $_{k[G]}^{k[G]}{\cal Y}{\cal D}$).
Note that $B\cap K=k1.$ Let
$\rho:B\raro K\ot B$ be the associated comodule structure and write
$\rho(x)=\sum x^1\ot x^2,$ $x\in B.$ 
By [R1], $B\cong A/AK^+$ as coalgebras, hence $B$ is a connected pointed 
coalgebra. Let $P(B):=\{x\in B|\Delta_B(x)=x\ot 1+1\ot x\}$ be the space 
of primitive elements of $B.$

\begin{lemma} \label{shtut1} 
For any $g\in G,$ $V_g=\{x\in P(B)|\rho(x)=g\ot x\}.$
\end{lemma}

\proof 
Let $x\in V_g.$ Since $g$ acts on $V_g$ by 
conjugation we may assume by [G1, Lemma 0.2], that $gx=\go xg$ for some
$1\ne \go\in k.$ Since $\pi(x)$ and $\pi(g)=g$ commute we 
must have that $\pi(x)=0.$ But then $(I\ot \pi)\Delta(x)=x\ot 1$ and 
hence $x\in B.$
Since $\Delta(x)=\sum x_1\times x_2^1\ot x_2^2\times 1,$ applying the
maps $\varepsilon\ot I\ot I\ot \varepsilon$ and 
$I\ot \varepsilon\ot I\ot \varepsilon$ to both sides of the equation 
$\sum x_1\times x_2^1\ot x_2^2\times 1=x\times 1\ot 1\times 1+1\times g\ot 
x\times 1,$ yields that $x\in P(B)$ and $\rho(x)=g\ot x$ as desired.

Suppose that $x\in P(B)$ satisfies $\rho(x)=g\ot x.$ Since
$\Delta(x)=x\ot 1+\rho(x),$ it follows that $x\in V_g$ as desired. \qed

\begin{lemma} \label{shtut2}
For every $x\in V_g,$ $x^2=0$ and $gx=-xg.$
\end{lemma}

\proof Suppose $V_g\ne 0$ and
let $0\ne x\in V_g.$ Then $S^2(x)=g^{-1}xg,$
$g^{-1}xg\ne x$ by [G1, Lemma 0.2], and $g^{-1}xg\in V_g.$ Since by 
Theorem \ref{ant}, $S^4=Id$ it follows that $g^2$ and $x$
commute, and
hence that $gx=-xg$ for every $x\in V_g.$

Second we wish to show that $x^2=0.$
By Lemma \ref{shtut1}, $x\in B$ and hence $x^2\in B$ ($B$ is a subalgebra
of $A$).
Since $\Delta(x^2)=x^2\ot 1+ g^2\ot x^2,$ and $x^2$ and
$g^2$ commute, it follows from [G1, Lemma 0.2] that 
$x^2=\ga(1-g^2)\in K$ for some $\ga\in k.$ We thus conclude that 
$x^2=0,$ as desired. \qed

Recall that the map $f_R:A^{*cop}\raro A$ is an isomorphism of Hopf
algebras, and let \linebreak
$F:G\times G\raro k^*$ be the associated non-degenerate skew
symmetric bilinear form on $G$ defined by
$F(g,h):=<f_R^{-1}(g),h>$ for every $g,h\in G.$

\begin{lemma} \label{shtut3}
For any $x\in V_g$ and $y\in V_h,$ $xy=F(h,g)yx.$  
\end{lemma}

\proof
If either $V_g=0$ or $V_h=0,$ there is nothing to prove. Suppose
$V_g,V_h\ne 0,$ and let $0\ne x\in V_g$ and $0\ne y\in V_h.$
Set $P:=f_R^{-1}(x).$ Then $P\in
P_{f_R^{-1}(g),\varepsilon}(A^{*cop}).$
Substituting $p:=P$ and $a:=y$ in equation (\ref{qt5}) yields that
$yx-F(g,h)xy=<P,y>(1-gh).$ Since $<P,y>(1-gh)\in B\cap K,$
it is equal to $0,$ and hence $yx=F(g,h)xy.$ \qed

\begin{lemma} \label{shtut4}
For any $a\in G$ and $x\in V_g,$ $xa=F(a,g)ax.$
\end{lemma}
 
\proof Set $P:=f_R^{-1}(x).$ Then the result follows by  
letting $p:=P$ and $a\in G$ in (\ref{qt5}), and noting that $P\in
P_{f_R^{-1}(g), \varepsilon}(A^{*cop}).$ \qed

We can now prove Theorem \ref{mainpoint}.

\noin
{\bf Proof of Theorem \ref{mainpoint}:}
Let $n:I_F\raro \Z^+$ be the nonnegative integer function
defined by
$n(g)=n_g,$ and let ${\cal D}:=(G,F,n).$ By [AEG, Theorem
6.1], $A$ is
generated as an algebra by $G\cup (\bigoplus_{g\in I_F}V_g).$
By Lemmas \ref{shtut2}-\ref{shtut4}, relations (\ref{3p}) are
satisfied.
Therefore there exists a surjection of Hopf algebras $\varphi:H({\cal
D})\raro A.$ Using Proposition 4.3 3), it is straightforward to verify
that $\varphi$ is injective on $P_{a,b}(H({\cal D}))$ for all
$a,b\in G.$ Since $H({\cal D})$ is pointed, $\varphi$ is also injective
(see e.g. [Mon, Corollary 5.4.7]). Hence $\varphi$ is an
isomorphism of
Hopf algebras. The rest of the theorem follows now from
Theorem \ref{quas}. \qed

\begin{remark}
{\rm Theorem 6.1 in [AEG] states
that a finite-dimensional cotriangular pointed Hopf
algebra is generated by its grouplike and skewprimitive
elements. This confirms the conjecture that this is the case
for {\em any} finite-dimensional pointed Hopf algebra over $\C$
[AS2], in the cotriangular case. The proof uses a categorical
point of view (alternatively the Lifting
method [AS1,AS2]).
}
\end{remark}

\section{Triangular Hopf algebras with the Chevalley property}

As we said in the introduction, semisimple cosemisimple
triangular Hopf algebras and minimal triangular pointed Hopf
algebras share in common the Chevalley property. 
In this section we describe the classification of finite-dimensional
triangular Hopf algebras with the Chevalley property, given in [AEG]. 

\subsection{Triangular Hopf algebras with Drinfeld
element of order $\le 2$}

We start by classifying triangular Hopf
algebras with $R-$matrix of rank $\le 2$.
We show that such a Hopf
algebra
is a suitable modification of a cocommutative Hopf
superalgebra
(i.e. the group algebra of a supergroup). On the other
hand, by Corollary \ref{kostantf},
a finite supergroup is a semidirect
product of
a finite group with an odd vector space on which this group
acts.

\subsubsection{The correspondence between Hopf algebras and
superalgebras}

We start with a correspondence theorem
between Hopf algebras and Hopf superalgebras.

\begin{ttheorem}\label{cores0}
There is a one to one correspondence between:
\ben
\item isomorphism classes of pairs $(A,u)$ where $A$ is an
ordinary Hopf
algebra, and $u$ is a grouplike element in $A$ such that
$u^2=1,$ and
\item
isomorphism classes of pairs $(\cH,g)$ where $\cH$
is a Hopf superalgebra, and
$g$ is a grouplike element in $\cH$
such that $g^{2}=1$ and $g x
g^{-1}=(-1)^{p(x)}x$ (i.e. $g$ acts on $x$ by
its parity),
\een
such that the tensor categories of representations of $A$ and
$\cH$ are equivalent.
\end{ttheorem}

\proof Let $(A,u)$ be an ordinary
Hopf algebra
with comultiplication $\Delta,$ counit $\eps$,
antipode $S$,
and a grouplike element $u$
such that $u^2=1.$
Let $\cH=A$ regarded as a superalgebra, where the
$\Z_2-$grading is
given by the adjoint
action of $u.$ For $a\in A,$ let us define $\Delta_0,
\Delta_1$ by writing $\Delta(a)=\Delta_0(a)+\Delta_1(a),$
where
$\Delta_0(a)\in A\ot A_0$ and $\Delta_1(a)\in A\ot A_1.$
Define a map
$\tilde\Delta:\cH\raro \cH\ot \cH$ by
$\tilde\Delta(a):=\Delta_0(a)-(-1)^{p(a)}(u\ot
1)\Delta_1(a).$ Define
$\tilde S(a):=u^{p(a)}S(a),$ $a\in A.$
Then it is straightforward
to verify that $(\cH,\tilde\Delta,\varepsilon,\tilde S)$
is a Hopf superalgebra.

The element $u$ remains grouplike in the new Hopf superalgebra,
and acts by parity, so we can set $g:=u$.

Conversely, suppose that $(\cH,g)$ is a pair where
$\cH$ is a Hopf superalgebra with
comultiplication $\tilde \Delta,$ counit $\varepsilon$,
antipode $\tilde S$, and a grouplike element $g$, with $g^2=1$,
acting by parity.
For $a\in \cH,$ let us define
$\tilde {\Delta}_0,\tilde {\Delta}_1$ by writing
$\tilde{\Delta}(a)=\tilde{\Delta}_0(a)+\tilde{\Delta}_1(a),$
where $\tilde{\Delta}_0(a)\in \cH\ot \cH_0$ and
$\tilde{\Delta}_1(a)\in \cH\ot \cH_1.$
Let $A=\cH$ as algebras, and
define a map
$\Delta:A\raro A\ot A$ by
$\Delta(a):=\tilde{\Delta}_0(a)-(-1)^{p(a)}(
g\ot 1)\tilde{\Delta}_1(a).$
Define $S(a):=g^{p(a)}\tilde S(a),$ $a\in A.$ Then it is
straightforward to verify that $(A,\Delta,\varepsilon,S)$ is an
ordinary Hopf
algebra, and we can set $u:=g$.

It is obvious that the two assignments constructed
above are inverse to each other.
The equivalence of tensor categories
is straightforward to verify. The theorem is proved.\qed

Theorem \ref{cores0} implies the following. Let $\cH$ be {\em
any} Hopf superalgebra, and $\C[\Z_2]\ltimes \cH$ be the
semidirect product, where the generator $g$ of $\Z_2$ acts on
$\cH$ by $gxg^{-1}=(-1)^{p(x)}x$. Then we can define an ordinary
Hopf algebra $\ol {\cH}$, which is the one corresponding to
$(\C[\Z_2]\ltimes \cH,g)$ under the correspondence of Theorem
\ref{cores0}.

The constructions of this subsection have the following explanation
in terms of Radford's biproduct construction [R1]. Namely
$\cH$ is
a Hopf algebra in the Yetter-Drinfeld category of $\C[\Z_2],$ so
Radford's biproduct construction yields a Hopf algebra structure
on $\C[\Z_2]\ot \cH,$ and it is straightforward to see that this
Hopf algebra is exactly $\ol {\cH}.$ Moreover, it is clear that
for any pair $(A,u)$ as in Theorem \ref{cores0}, $gu$ is central
in $\ol {\cH}$ and $A=\ol {\cH}/(gu-1).$

\subsubsection{Correspondence of twists}

Let us say that a twist $J$ for a Hopf algebra $A$
with an involutive grouplike element
$g$ is {\it even} if
it is invariant under $\Ad(g)$.

\begin{pproposition}\label{1a}
Let $(\cH,g)$ be a pair as in Theorem
\ref{cores0}, and let
$A$ be the associated ordinary Hopf algebra. Let
${\cal J}\in \cH\ot \cH$ be an even element. Write
${\cal J}=
{\cal J}_0+{\cal J}_1,$ where ${\cal J}_0\in \cH_0\ot
\cH_0$ and
${\cal J}_1\in \cH_1\ot \cH_1.$ Define
$J:={\cal J}_0-
(g\ot 1) {\cal J}_1$. Then $J$ is an even twist for $A$
if and only if ${\cal J}$ is a twist for $\cH.$
Moreover, $\cH^{{\cal J}}$ corresponds to $A^J$ under the
correspondence in Theorem \ref{cores0}.
Thus, there is a one to one correspondence between even twists
for $A$ and twists for $\cH$, given by $J\to {\cal J}$.
\end{pproposition}

\proof Straightforward.\qed

\subsubsection{The Correspondence between triangular Hopf
algebras and superalgebras}

Let us now return to our main subject, which is
triangular Hopf algebras and superalgebras.
For triangular Hopf algebras whose Drinfeld element $u$ is
involutive, we will make the natural choice of the
element $u$ in Theorem \ref{cores0}, namely define it to be the
Drinfeld element of $A$.

\begin{ttheorem}\label{cores}
The correspondence of Theorem \ref{cores0}
extends to a one to one correspondence between:
\ben
\item isomorphism classes of ordinary
triangular Hopf
algebras $A$ with Drinfeld element $u$ such that $u^2=1,$ and
\item
isomorphism classes of pairs $(\cH,g)$ where $\cH$
is a
triangular Hopf superalgebra with Drinfeld element $1$ and
$g$ is an element of ${\bf G}(\cH)$
such that $g^{2}=1$ and $gxg^{-1}=(-1)^{p(x)}x.$
\een
\end{ttheorem}

\proof Let $(A,R)$ be a triangular Hopf algebra with $u^2=1$.
Since $(S\ot S)(R)=R$ and $S^2=\Ad(u)$ [Dr2], $u\ot u$ and $R$
commute. Hence we can write $R=R_0+R_1$, where $R_0\in A_0\ot
A_0$ and $R_1\in A_1\ot A_1.$ Let $\cR:=(R_0+(1\ot u)R_1)R_u$.
Then $\cR$ is even. Indeed, since $R_0=1/2(R+(u\ot 1)R(u\ot 1))$
and $R_1=1/2(R-(u\ot 1)R(u\ot 1)),$ $u\ot u$ and $\cR$ commute.

It is now straightforward to show that $(\cH,\cR)$ is triangular
with Drinfeld element $1.$ Let us show for instance that $\cR^{-1}=
\cR_{21}.$ Let us use the
notation $a*b, X^{21}$ for multiplication
and opposition in the tensor square of a superalgebra, and the
notation $ab,X^{op}$ for usual algebras. Then,
$$\cR*\cR_{21}=
(R_0+(1\ot u)R_1)R_u*(R_0^{op} - (u\ot 1)R_1^{op})R_u.$$ Since,
$R_uR_0=R_0R_u,$ $R_uR_1=-(u\ot u)R_1R_u,$ we get that the RHS
equals
$$(R_0+(1\ot u)R_1)*(R_0^{op}+(1\ot u)R_1^{op})=
R_0R_0^{op}+R_1R_1^{op}+(1\ot u)(R_1R_0^{op}+R_0R_1^{op}).$$
But, $R_0R_0^{op}+R_1R_1^{op}=1$ and $(1\ot
u)(R_1R_0^{op}+R_0R_1^{op})=0,$ since $RR^{op}=1,$ so we
are done.

Conversely, suppose that $(\cH,g)$ is a pair where
$\cH$ is a triangular Hopf superalgebra with $R-$matrix
$\cR$
and Drinfeld element $1.$ Let
$\cR=\cR_0+\cR_1$, where
$\cR_0$ has even components, and $\cR_1$ has odd
components. Let $R:=(\cR_0+(1\ot g)\cR_1) R_g$.
Then it is straightforward
to show that $(A,R)$ is
triangular with Drinfeld element $u=g$.
The theorem is proved.\qed

\begin{ccorollary}\label{mintr}
If $(\cH,\cR)$ is a triangular
Hopf superalgebra with Drinfeld element 1,
then the Hopf algebra $\overline{\cH}$ is also
triangular, with the $R-$matrix
\begin{equation}\label{olr}
\ol R:=(\cR_0+(1\otimes g)\cR_1)R_g,
\end{equation}
where $g$ is the grouplike element adjoined to $\cH$ to
obtain $\overline{\cH}$. Moreover, $\cH$ is minimal if and
only if so is $\ol {\cH}.$
\end{ccorollary}

\proof Clear.\qed

The following corollary, combined with Kostant's theorem, gives a 
classification of triangular Hopf algebras
with $R$-matrix of rank $\le 2$ (i.e. of the form
$R_u$ as in (\ref{ru}), where $u$ is a grouplike of order
$\le 2$).

\begin{ccorollary}\label{cores'}
The correspondence of Theorem \ref{cores}
restricts to a one to one correspondence between:
\ben
\item isomorphism classes of ordinary
triangular Hopf algebras with $R-$matrix of
rank $\le 2,$ and
\item
isomorphism classes of pairs $(\cH,g)$ where
$\cH$ is a cocommutative Hopf superalgebra and $g$ is an
element of ${\bf G}(\cH)$
such that $g^{2}=1$ and $g x
g^{-1}=(-1)^{p(x)}x$.
\een
\end{ccorollary}
\proof Let $(A,R)$ be an ordinary triangular Hopf algebra with
$R-$ matrix of rank $\le 2.$ In particular, the Drinfeld element
$u$ of $A$ satisfies $u^2=1,$ and $R=R_u$. Hence by Theorem
\ref{cores}, $(\cH,\tilde\Delta,\cR)$ is a triangular Hopf
superalgebra. Moreover, it is cocommutative since $\cR=R_uR_u=1.$

Conversely, for any $(\cH,g)$,
by Theorem \ref{cores}, the pair $(A,R_{g})$ is an
ordinary triangular Hopf algebra, and clearly the
rank of $R_{g}$ is $\le 2.$ \qed

In particular, Corollaries \ref{kostantf} and \ref{cores'}
imply that finite-dimensional
triangular Hopf algebras with $R-$matrix of rank $\le 2$
correspond to supergroup algebras.
In view of this, we make the following definition.

\begin{ddefinition}
A finite-dimensional triangular Hopf algebra with $R-$matrix of
rank $\le 2$ is called a modified supergroup algebra.
\end{ddefinition}

\subsubsection{Construction of twists for supergroup
algebras}

\begin{pproposition}\label{st}
Let $\cH=\C[G] \ltimes \Lambda V$ be a supergroup algebra.
Let $r\in S^2V.$ Then ${\cal J}:=e^{r/2}$ is a twist for $\cH.$
Moreover,
$((\Lambda V)^{{\cal J}},{\cal J}_{21}^{-1}{\cal J})$ is minimal
triangular if
and only if $r$ is nondegenerate.
\end{pproposition}
\proof Straightforward. \qed

\begin{eexample}
{\rm Let $G$ be the group of order $2$ with generator $g.$
Let $V:=\C$ be the nontrivial $1-$dimensional
representation of $G,$ and write
$\Lambda V=sp\{1,x\}.$ Then the associated ordinary
triangular Hopf algebra
to $(\cH,g):=(\C[G] \ltimes \Lambda V,g)$ is
Sweedler's $4-$dimensional
Hopf algebra $A$ [Sw] (see Example \ref{sweed}) with the triangular
structure
$R_g.$ It is known [R2] that
the set of triangular structures on $A$ is parameterized
by $\C;$ namely,
$R$ is a triangular structure on $A$ if and only if
$$R=R_{\lambda}:=
R_g-\frac{\lambda}{2}(x\ot x-gx\ot x+x\ot gx +gx\ot gx),
\;\lambda\in \C.$$ Clearly,
$(A,R_{\lambda})$ is minimal if and only if $\lambda\ne 0.$

Let $r\in S^2V$ be defined by $r:=\lambda x\otimes x,$
$\lambda\in \C.$ Set
${\cal J}_\lambda:=e^{r/2}=1+\frac{1}{2}\lambda x\ot x$;
it is a
twist for $\cH.$
Hence, $J_{\lambda}:=1-\frac{1}{2}\lambda gx\ot x$ is a
twist for $A.$ It is
easy to check that
$R_{\lambda}=(J_{\lambda})_{21}^{-1}R_gJ_{\lambda}.$
Thus, $(A,R_\lambda)=(A,R_0)^{J_\lambda}$.
}
\end{eexample}

\begin{rremark}{\rm  In fact, Radford's classification
of triangular structures on $A$ can be easily deduced
from Lemma \ref{5c} below.
}
\end{rremark}

\subsection{The Chevalley property}

Recall that in the introduction we made the following
definition.

\begin{definition} A Hopf algebra $A$ over $\C$ is said to have
the Chevalley property if the tensor product of any two simple
$A$-modules is semisimple. More generally, let us say that a
tensor category has the Chevalley property if the tensor product
of two simple objects is semisimple.
\end{definition}

Let us give some equivalent formulations of the Chevalley
property.

\begin{proposition}\label{chev}
Let $A$ be a finite-dimensional Hopf algebra over $\C.$ The following
conditions  are equivalent:
\ben
\item $A$ has the Chevalley property.
\item The category of (right) $A^*$-comodules has the
Chevalley property.
\item $\Corad(A^*)$ is a Hopf subalgebra of $A^*.$
\item $\Rad(A)$ is a Hopf ideal and thus $A/\Rad(A)$
is a Hopf algebra.
\item $S^2=Id$ on $A/\Rad(A)$, or equivalently on $\Corad(A^*)$.
\een
\end{proposition}

\proof
(1. $\Leftrightarrow$ 2.) Clear, since the categories of left
$A$-modules
and right $A^*$-comodules are equivalent.

(2. $\Rightarrow$ 3.)
Recall the definition of a matrix coefficient of a comodule $V$
over $A^*$. If $\rho: V \to V \otimes A^*$ is the coaction, $v\in
V$, $\alpha \in V^*$, then
$$
\phi^V_{v, \alpha}:=(\alpha \otimes Id)\rho (v) \in A^*.
$$
It is well-known that:

(a) The coradical of $A^*$ is the linear span of the matrix
coefficients of
all simple $A^*$-comodules.

(b) The product in $A^*$ of two matrix coefficients is a matrix
coefficient of the tensor product. Specifically,
$$
\phi^V_{v, \alpha} \phi^W_{w, \beta} = \phi^{V\otimes
W}_{v\otimes w,
\alpha\otimes \beta}.
$$
It follows at once from (a) and (b) that $\Corad(A^*)$ is a
subalgebra of $A^*.$
Since the coradical is stable under the antipode, the claim
follows.

(3. $\Leftrightarrow$ 4.) To say that
$\Rad(A)$ is a Hopf ideal is equivalent to
saying that $\Corad(A^*)$ is a Hopf algebra, since
$\Corad(A^*)=(A/\Rad(A))^*.$

(4. $\Rightarrow$ 1.) If $V,W$ are simple
$A-$modules then they factor
through
$A/\Rad(A).$ But $A/\Rad(A)$ is a Hopf algebra,
so $V\otimes W$ also factors through
$A/\Rad(A),$ so it is semisimple.

(3. $\Rightarrow$ 5.) Clear, since a cosemisimple Hopf
algebra is involutory.

(5. $\Rightarrow$ 3.) Consider the subalgebra $B$ of $A^*$
generated by $\Corad(A^*).$
This is a Hopf algebra, and $S^2=Id$ on it.
Thus, $B$ is cosemisimple and hence
$B=\Corad(A^*)$ is a Hopf subalgebra of $A^*$.
\qed

\begin{remark} {\rm The assumption that the base field has
characteristic $0$
is needed only in the proof of (5. $\Leftrightarrow$ 3.)
}
\end{remark}

\subsection{The classification of triangular Hopf algebras
with the Chevalley property}

\subsubsection{The main theorem}

The main result of Section 7 is the following theorem.

\begin{ttheorem}\label{mainchev} Let $A$ be a
finite-dimensional triangular
Hopf algebra over $\C.$ Then the following are equivalent:
\ben
\item $A$ is a twist of a finite-dimensional triangular
Hopf algebra with $R-$matrix of
rank $\le 2$ (i.e. of a modified supergroup algebra).

\item $A$ has the Chevalley property.
\een
\end{ttheorem}

The proof of this theorem is contained in the next
subsubsection.

\subsubsection{Proof of the main theorem}

We shall need the following result whose proof is given in [AEG].

\begin{ttheorem}\label{1} Let $\cH$ be a local
finite-dimensional Hopf
superalgebra (not necessarily supercommutative).
Then $\cH=\Lambda V^*$ for a finite-dimensional vector
space $V.$
In other words, $\cH$ is the function algebra of an odd
vector space $V$.
\end{ttheorem}

\begin{rremark} {\rm Note that in the commutative case
Theorem \ref{1}
is a special case of Proposition 3.2 of [Ko].
}
\end{rremark}

We start by giving a super-analogue of Theorem 3.1 in [G4].

\begin{llemma}\label{super3.1}
Let $\cH$ be a minimal triangular pointed Hopf
superalgebra.
Then $\Rad(\cH)$ is a Hopf ideal, and $\cH/\Rad(\cH)$ is
minimal triangular.
\end{llemma}

\proof The proof is a tautological generalization of the
proof of Theorem 3.1 in [G4] to the super case.

First of all, it is clear that $\text{Rad}(\cH)$ is a Hopf
ideal,
since its orthogonal complement (the coradical of $\cH^*$)
is a sub Hopf superalgebra (as $\cH^*$ is isomorphic to
$\cH^{cop}$
as a coalgebra, and hence is pointed).
Thus, it remains to show that the triangular structure on
$\cH$
descends to a minimal triangular structure on
$\cH/\Rad(\cH)$.
For this, it suffices to prove that the
composition of the Hopf superalgebra maps
$$
\text{Corad}(\cH^{*cop})\hookrightarrow \cH^{*cop}\to \cH\to
\cH/\Rad(\cH)
$$
(where the middle map is given by the $R-$matrix) is an
isomorphism. But this follows from the fact that
for any surjective coalgebra map $\eta: C_1\to C_2$,
the image of the coradical of $C_1$ contains the coradical of
$C_2$ (see e.g. [Mon, Corollary 5.3.5]): One needs to apply
this
statement to the map $\cH^{*cop}\to \cH/\text{Rad}(\cH)$.
\qed

\begin{llemma}\label{3} Let $\cH$ be a minimal
triangular pointed Hopf superalgebra, such that the \linebreak
$R-$matrix $\cR$ of $\cH$ is unipotent (i.e. $\cR-1\ot 1$ is $0$
in $\cH/\Rad(\cH)\ot \cH/\Rad(\cH)$). Then $\cH= \Lambda V$ as a
Hopf superalgebra, and $\cR=e^r,$ where $r\in S^2V$ is a
nondegenerate symmetric (in the usual sense) bilinear form on
$V^*$.
\end{llemma}

\proof By Lemma \ref{super3.1}, $\Rad(\cH)$ is a Hopf ideal, and
$\cH/\Rad(\cH)$ is minimal triangular. But the $R-$matrix of
$\cH/\Rad(\cH)$ must be $1\ot 1,$ so $\cH/\Rad(\cH)$ is
$1-$dimensional. Hence $\cH$ is local, so by Theorem \ref{1},
$\cH=\Lambda V.$ If $\cR$ is a triangular structure on $\cH$ then
it comes from an isomorphism $\Lambda V^*\raro \Lambda V$ of Hopf
superalgebras, which is induced by a linear isomorphism
$r:V^*\raro V.$ So $\cR=e^r,$ where $r$ is regarded as an element
of $V\ot V.$ Since $\cR \cR_{21}=1,$ we have $r+r^{21}=0$ (where
$r^{21}=-r^{op}$ is the opposite of $r$ in the supersense), so
$r\in S^2V$. \qed

\begin{rremark}
{\rm The classification of pointed finite-dimensional Hopf
algebras with coradical of dimension $2$ is known [CD,N]. In [AEG] we
used the Lifting method [AS1,AS2] to give an
alternative proof. Below we shall need the following more precise
version of this result in the triangular case. }
\end{rremark}

\begin{llemma}\label{5c} Let $A$ be a
minimal triangular pointed Hopf algebra,
whose coradical is $\C[\Z_2]=sp\{1,u\},$ where $u$ is the
Drinfeld element of $A.$
Then $A=\overline{(\Lambda V)^{{\cal J}}}$
with the triangular structure of Corollary \ref{mintr},
where
${\cal J}=e^{r/2},$ with $r\in S^2V$ a nondegenerate element.
In particular, $A$ is a twist of a modified supergroup algebra.
\end{llemma}

\proof Let $\cH$ be the associated triangular Hopf superalgebra
to $A,$ as described in Theorem \ref{cores}. Then the
$R-$matrix of $\cH$
is unipotent, because it turns into
$1\ot 1$ after killing the radical.

Let ${\cH}_m$ be the minimal part of $\cH.$ By
Lemma \ref{3}, ${\cH}_m=\Lambda V$ and $\cR=e^r,$ $r\in
S^2V.$ So if ${\cal J}:=e^{r/2}$ then
$\cH^{{\cal J}^{-1}}$ has $R-$matrix equal to
$1\ot 1.$ Thus, $\cH^{{\cal J}^{-1}}$ is cocommutative,
so by Corollary \ref{kostantf}, it equals
$\C[\Z_2]\ltimes \Lambda V.$ Hence
$\cH=\C[\Z_2]\ltimes (\Lambda V)^{{\cal J}}$,
and the result follows from Proposition \ref{1a}.
\qed

We shall need the following lemma.

\begin{llemma}\label{jh}
Let $B\subseteq A$ be finite-dimensional associative unital
algebras.
Then any simple $B-$module is a constituent
(in the Jordan-Holder series) of some simple $A-$module.
\end{llemma}

\proof Since $A,$ considered as a $B-$module, contains $B$
as a $B-$module, any simple $B-$module is a constituent of
$A.$

Decompose $A$ (in the Grothendieck group of $A$) into
simple $A-$modules: $A=\sum V_i.$ Further decomposing as
$B-$modules, we get $V_i=\sum W_{ij},$ and hence
$A=\sum_i\sum_j W_{ij}.$ Now, by Jordan-Holder theorem,
since $A$ (as a $B-$module) contains all simple
$B-$modules, any simple $B-$module $X$ is in $\{W_{ij}\}.$
Thus, $X$ is a constituent of some $V_i,$ as desired. \qed

\begin{pproposition}\label{6} Any minimal
triangular Hopf algebra $A$ with the Chevalley property
is a twist of a triangular Hopf algebra with
$R-$matrix of rank $\le 2.$
\end{pproposition}

\proof
By Proposition \ref{chev}, the coradical $A_0$ of $A$ is a Hopf
subalgebra,
since $A \cong A^{*cop}$, being minimal triangular.
Consider the  Hopf algebra map $\varphi: A_{0} \to A^{*cop}/
\Rad(A^{*cop})$,
given by the composition of the following maps:
$$
A_{0} \hookrightarrow A \cong A^{*cop} \to A^{*cop}/
\Rad(A^{*cop}),
$$
where the second map is given by the $R$-matrix. We claim that
$\varphi$ is an
isomorphism. Indeed, $A_{0}$ and $A^{*cop}/ \Rad(A^{*cop})$ have the same
dimension, since
$\Rad(A^{*cop}) = (A_{0})^{\bot}$, and $\varphi$ is injective,
since
$A_{0}$
is semisimple by
[LR]. Let $\pi: A\to A_0$ be the associated projection.

We see, arguing exactly as in [G4, Theorem 3.1], that  $A_0$ is
also minimal triangular, say with $R-$matrix $R_0$.

Now, by [EG1, Theorem 2.1], we can find a twist $J$ in $A_0\ot A_0$
such that $(A_0)^{J}$ is isomorphic to a group algebra
and has $R-$matrix $(R_0)^{J}$ of rank $\le
2$. Notice that here we are relying on Deligne's theorem, as mentioned in the
introduction.

Let us now consider $J$ as an element of $A_0\ot A_0$ and the twisted
Hopf algebra
$A^{J}$, which is again triangular.

The projection $\pi:A^{J} \to (A_0)^{J}$ is still a Hopf algebra map, and
sends
$R^J$ to $(R_0)^J.$ It induces a projection $(A^J)_m\to \C[\Z_2],$ whose
kernel
$K_m$ is contained in the kernel of $\pi$.
Because any simple $(A^J)_m-$module
is contained as a constituent in a simple $A-$module (see
Lemma \ref{jh}), $K_m=\Rad((A^J)_m).$ Hence,
$(A^J)_m$ is minimal triangular and $(A^J)_m/\Rad((A^J)_m)=(\C[\Z_2],
R_u).$ It follows, again by minimality, that $(A^J)_m$ is also pointed
with coradical
isomorphic to $\C[\Z_2]$.
So by Lemma \ref{5c}, $(A^J)_m,$ and hence $A^J,$ can be
further twisted
into a triangular Hopf algebra with $R-$matrix of rank
$\le 2,$ as desired. \qed

Now we can prove the main theorem.

\noin
{\bf Proof of theorem \ref{mainchev}:}
(2. $\Rightarrow$ 1.) By Proposition
\ref{chev},
$A/\Rad(A)$ is a semisimple Hopf algebra. Let
$A_m$ be the minimal part
of $A,$ and $A_m'$ be the image of $A_m$ in
$A/\Rad(A).$ Then $A_m'$ is
a semisimple Hopf algebra.

Consider the kernel $K$ of the projection $A_m\to
A_m'.$
Then $K=\Rad(A)\cap A_m.$ This means that any
element $k\in K$ is
zero in any simple $A-$module. This implies that
$k$ acts by zero in any simple $A_m-$module,
since by Lemma \ref{jh}, any simple $A_m-$module
occurs as a constituent of some simple $A-$module. Thus,
$K$ is contained in $\Rad(A_m).$
On the other hand, $A_m/K$ is semisimple, so
$K=\Rad(A_m).$ This shows that $\Rad(A_m)$ is a Hopf
ideal. Thus, $A_m$ is
minimal triangular
satisfying the conditions of Proposition \ref{6}.
By Proposition \ref{6}, $A_m$ is a twist of a triangular Hopf algebra
with $R-$matrix of rank $\le 2$. Hence $A$ is a twist of a
triangular Hopf algebra with $R-$matrix of rank $\le 2$
(by the same twist), as desired.

(1. $\Rightarrow$ 2.) By assumption, $\Rep(A)$ is
equivalent to $\Rep(\tilde G)$ for some supergroup $\tilde G$
(as a tensor category without braiding).
But we know that supergroup algebras have the Chevalley
property, since, modulo their radicals, they are group
algebras.  This concludes the proof of the main
theorem. \qed

\begin{rremark} {\rm Notice that it follows from the proof
of the
main
theorem
that any triangular Hopf algebra with the Chevalley
property can be
obtained by twisting of a triangular Hopf algebra
with $R-$matrix of rank $\le 2$ by an {\it even} twist.
}
\end{rremark}

\begin{ddefinition} If a triangular Hopf algebra
$A$ over $\C$ satisfies condition 1. or 2. of Theorem
\ref{mainchev}, we will
say that $H$ is of supergroup type.
\end{ddefinition}

\subsubsection{Corollaries of the main theorem}

\begin{ccorollary} \label{minim} A finite-dimensional triangular
Hopf
algebra $A$ is of supergroup type if and only if so is
its minimal part $A_m.$
\end{ccorollary}

\proof If $A$ is of supergroup type then
$\Rad(A)$ is a Hopf
ideal, so like
in the proof of Theorem \ref{mainchev} (2. $\Rightarrow$ 1.) we
conclude
that $\Rad(A_m)$
is
a Hopf ideal, i.e. $A_m$ is of supergroup type.

Conversely, if $A_m$
is of supergroup type
then $A_m$ is a twist of a triangular Hopf algebra
with $R-$matrix of rank $\le 2$. Hence $A$ is a twist of a
triangular
Hopf algebra with $R-$matrix of rank $\le 2$
(by the same twist), so $A$ is
of supergroup type. \qed

\begin{ccorollary}\label{6a} A finite-dimensional
triangular Hopf algebra
whose coradical is a Hopf subalgebra is of supergroup type.
In particular, this is the case for finite-dimensional triangular 
pointed Hopf algebras.
\end{ccorollary}

\proof This follows from Corollary \ref{minim}.\qed

\begin{ccorollary} \label{triangcopoint} Any
finite-dimensional triangular
basic Hopf algebra is of supergroup type.
\end{ccorollary}

\proof A basic Hopf algebra automatically has
the Chevalley property
since all its irreducible modules are
$1-$dimensional. Hence the result follows from the
main theorem. \qed

\subsection{Categorical dimensions in symmetric
categories with finitely many irreducibles are integers}

In [AEG] we classified finite-dimensional
triangular Hopf algebras
with the Chevalley property.
We also gave one result that is valid in
a greater generality for any finite-dimensional
triangular Hopf
algebra, and even for any symmetric rigid category
with finitely many irreducible objects.

Let ${\cal C}$ be a $\C-$linear abelian symmetric rigid
category with ${\bold 1}$ as its unit
object, and suppose that $\End({\bold 1})=\C.$ Recall
that there is a natural
notion of dimension in ${\cal C},$ generalizing the
ordinary dimension of an object
in $\Vect,$ and having the properties of being additive and
multiplicative with respect to the tensor product. 
Let $\beta$ denote the
commutativity constraint in ${\cal
C},$ and for an object
$V,$ let $ev_V,$ $coev_V$ denote the associated
evaluation and coevaluation
morphisms.

\begin{definition}\label{dimcat} {\bf [DM]}
The categorical dimension $\text
{dim}_c(V)\in \C$ of $V\in {\cal C}$ is the
morphism
\begin{equation}\label{dimc}
\text {dim}_c(V):{\bold 1}\stackrel{ev_V}{\longrightarrow}V\ot
V^*\stackrel{\beta_{V,V^*}}{\longrightarrow}V^*\ot
V\stackrel{coev_V}{\longrightarrow}{{\bold 1}}.
\end{equation}
\end{definition}

The main result of this subsection is the following:

\begin{theorem}\label{cd} In any $\C-$linear abelian symmetric
rigid tensor category ${\cal C}$
with finitely
many irreducible objects, the categorical dimensions of
objects are
integers.
\end{theorem}

\proof First note that the categorical dimension of any
object $V$ of ${\cal C}$ is an
algebraic integer. Indeed, let $V_1\dots,V_n$ be
the irreducible objects of ${\cal
C}.$ Then $\{V_1\dots,V_n\}$ is a basis of the
Grothendieck ring of ${\cal C}.$
Write $V\otimes V_i=\sum_j N_{ij}(V)V_j$ in the
Grothendieck ring. Then $N_{ij}(V)$ is
a matrix with integer entries, and $\dim_c(V)$ is
an eigenvalue of this matrix.
Thus, $\dim_c(V)$ is an algebraic integer.

Now, if $\dim_c(V)=d$ then it is easy to show (see e.g.
[De1]) that
$$\dim_c(S^kV)=d(d+1)\cdots (d+k-1)/k!,$$
and
$$\dim_c(\Lambda^k V)=d(d-1)\cdots (d-k+1)/k!,$$
hence they are also algebraic
integers. So the theorem follows from:

\noin
{\bf Lemma.} {\em Suppose $d$ is an algebraic
integer such that $d(d+1)\cdots (d+k-1)/k!$
and \linebreak $d(d-1)\cdots (d-k+1)/k!$ are algebraic
integers for all $k.$ Then $d$ is an
integer.}

\noin
{\bf Proof:} Let $Q$ be the minimal monic polynomial of
$d$ over $\Z.$
Since $d(d-1)\cdots (d-k+1)/k!$ is an algebraic
integer, so are the
numbers $d'(d'-1)\cdots (d'-k+1)/k!,$ where $d'$ is
any algebraic conjugate of $d$.
Taking the product over all conjugates, we get that
$$N(d)N(d-1)\cdots N(d-k+1)/(k!)^n$$
is an integer,
where $n$ is the degree of
$Q.$
But $N(d-x)=(-1)^nQ(x).$ So we get that
$Q(0)Q(1)\cdots Q(k-1)/(k!)^n$ is an integer.
Similarly from the identity for $S^kV,$ it
follows that $Q(0)Q(-1)\cdots Q(1-k)/(k!)^n$ is an
integer. Now, without loss of generality, we can
assume that
$Q(x)=x^n+ax^{n-1}+...,$ where $a\le 0$
(otherwise replace $Q(x)$ by
$Q(-x);$ we can do it since our condition is
symmetric under this
change).
Then for large $k,$ we have $Q(k-1)<k^n,$ so the
sequence
$b_k:=Q(0)Q(1)\cdots Q(k-1)/k!^n$ is decreasing in
absolute value
or zero starting from some place. But a sequence
of integers cannot
be strictly decreasing in absolute value forever.
So $b_k=0$ for some
$k,$ hence $Q$ has an integer root. This means
that $d$ is an integer
(i.e. $Q$ is linear), since $Q$ must be irreducible
over the rationals.
This concludes the proof of the lemma, and hence
of the theorem. \qed

\begin{corollary}\label{cdtr} For any triangular Hopf
algebra $A$ (not
necessarily
finite-dimensional), the categorical dimensions of
its finite-dimensional
representations are integers.
\end{corollary}

\proof It is enough to consider the minimal part $A_m$ of
$A$ which is
finite-dimensional, since $\dim_c(V)=\tr(u_{|V})$ for any
module $V$ (where $u$ is the Drinfeld element of $A$), and
$u\in A_m.$
Hence the result follows from Theorem \ref{cd}.
\qed

\begin{remark} {\rm Theorem \ref{cd} is false without the
finiteness conditions. In fact, in this case any complex
number can be a dimension, as is
demonstrated in examples constructed by Deligne [De2,
p.324-325]. Also, it is well known that the theorem
is false for ribbon, nonsymmetric categories
(e.g. for fusion categories of
semisimple representations of finite-dimensional quantum
groups at roots of unity [L],
where dimensions can be irrational algebraic integers).
}
\end{remark}

\begin{remark} {\rm In any rigid braided tensor category with
finitely many
irreducible objects, one can define the Frobenius-Perron
dimension
of an object $V$, $\FPdim(V),$ to be the largest positive
eigenvalue of the
matrix of multiplication by $V$ in the Grothendieck ring.
This dimension is well defined
by the Frobenius-Perron theorem, and has the usual
additivity and
multiplicativity properties.
For example, for the category of representations of a
quasi-Hopf algebra, it is just the usual dimension of the
underlying vector space.
If the answer to
Question \ref{q2a} is positive then $\FPdim(V)$
for symmetric categories is always an integer,
which is equal to $\dim_c(V)$ modulo $2.$ It would be
interesting
to check this, at least in the case of modules over a
triangular
Hopf algebras, when the integrality of $\FPdim$ is automatic
(so only the mod $2$ congruence has to be checked).}
\end{remark}

\section{Questions}

We conclude the paper with some natural questions motivated
by the above results [AEG,G4].

\begin{Question}\label{q}
{\rm 
Let $(A,R)$ be {\em any} finite-dimensional triangular
Hopf algebra with Drinfeld element $u.$ Is it true that
$S^4=Id$? Does $u$ satisfy $u^2=1$? Is it true that
$S^4=Id$ implies $u^2=1$?
}
\end{Question}

\begin{Remark} {\rm
A positive answer to the second question in Question
\ref{q} will imply that an
odd-dimensional triangular Hopf algebra must be semisimple.

Note that if $A$ is of supergroup type, then the
answer to Question \ref{q} is positive. Indeed,
since $S^2=Id$ on the semisimple 
part of $A,$ $u$ acts by a scalar in any irreducible
representation of
$A.$ In fact, since $\tr(u)=\tr(u^{-1}),$ we have that $u=1$
or $u=-1$ on
any irreducible representation of $A,$ and hence $u^2=1$ on any
irreducible
representation of $A.$ Thus, $u^2$ is unipotent. But it is of
finite order (as it is a grouplike element), so it is equal
to $1$ as desired.
}
\end{Remark}

\begin{Question} \label{q1}
Does any finite-dimensional
triangular Hopf algebra over $\C$ have the Chevalley property
(i.e. is of supergroup type)?
Is it true under the assumption that $S^4=Id$
or at least under the assumption that $u^2=1$?
\end{Question}

\begin{Remark} {\rm Note that the answer to question \ref{q1}
is negative in the infinite dimensional case. Namely, although the
answer is
positive
in the cocommutative case (by [C]), it is negative already
for triangular Hopf algebras with $R-$matrix of rank $2,$
which
correspond to cocommutative Hopf superalgebras.
Indeed, let us take the cocommutative Hopf superalgebra
$\cH:=\text{U}(\text{gl}(n|n))$
(for the definition of the Lie superalgebra
$\text{gl}(n|n),$ see [KaV, p.29]). The associated
triangular Hopf algebra $\ol {\cH}$ does not have the
Chevalley property,
since it is well known that Chevalley theorem fails
for Lie superalgebras (e.g. $\text{gl}(n|n)$); more
precisely, already
the product of the vector and covector representations
for this Lie superalgebra is not semisimple.
}
\end{Remark}

\begin{Remark} {\rm It follows from Corollary \ref{minim}
that a positive answer to Question \ref{q1} in the minimal case
would imply the general positive answer.
}
\end{Remark}

Here is a generalization of Question \ref{q1}.

\begin{Question}\label{q2} {\rm Does any $\C-$linear abelian
symmetric rigid tensor category, with $\text{End}({\bf 1})=\C$
and finitely many
simple objects, have the Chevalley property? }
\end{Question}

Even a more ambitious question:

\begin{Question}\label{q2a}{\rm
Is such a category equivalent to the category of
representations of a
finite-dimensional triangular
Hopf algebra with $R-$matrix of rank $\le 2$? In particular,
is it equivalent to the category of representations of a
supergroup, as a category without braiding?
Are these statements valid at least for categories
with Chevalley property?
For semisimple categories?}
\end{Question}

\begin{Remark} {\rm Note that Theorem \ref{cd}
can be regarded as a piece of supporting evidence for a
positive answer to Question \ref{q2a}.
}
\end{Remark}


\begin{thebibliography}{[EGGS]} 
\bibitem
[AEG]{aeg} N. Andruskiewitsch, P. Etingof and S. Gelaki,
Triangular Hopf 
Algebras with the Chevalley Property, {\em submitted},
math.QA/0008232.
\bibitem
[AS1]{as1} N. Andruskiewitsch and H.-J. Schneider, Lifting
of Quantum
Linear Spaces and Pointed Hopf Algebras of Order $p^3$, {\em J.
Algebra} {\bf 209} (1998), 658-691.
\bibitem
[AS2]{AS2}
N. Andruskiewitsch and H.-J. Schneider, Finite quantum
groups and Cartan matrices, {\em Adv. in  Math.} {\bf 154} (2000),
1-45.
\bibitem
[AT]{at} J.H. Conway, R.T. Curtis, S.P. Norton,
R.A. Parker and R.A. Wilson, Atlas of Finite Groups,
{\em Clarendon Press, Oxford} (1985).
\bibitem
[C]{c} C. Chevalley, Theory of Lie groups, v.III,
1951 (in French).
\bibitem
[CD]{CD} S. Caenepeel and  S. D\u{a}sc\u{a}lescu,
On pointed Hopf algebras of dimension $2\sp n$,
{\em Bull. London Math. Soc.} {\bf 31} (1999),
17-24
\bibitem
[CR]{cr} C. Curtis and I. Reiner, Methods of
Representation Theory {\bf
1}, John Wiley \& Sons,
Inc. (1981).
\bibitem
[De1]{de1} P. Deligne, Categories Tannakiennes, In The
Grothendick
Festschrift, Vol. II, Prog. Math. {\bf 87} (1990), 111-195.
\bibitem
[De2]{de2} P. Deligne, La s\'erie exceptionnelle de
groupes de Lie.
(French) [The exceptional series of Lie groups],
{\em C. R. Acad. Sci.
Paris Sor. I Math.} {\bf 322} (1996), no.4, 321--326.
\bibitem 
[DM]{dm} P. Deligne and J. Milne, Tannakian
Categories, Lecture
Notes in Mathematics {\bf 900}, 101-228, 1982.
\bibitem
[Dr1]{dr1} V. Drinfeld, Quantum Groups, {\em
Proceedings of the
International Congress of Mathematics, Berkeley}
(1987), 798-820.
\bibitem
[Dr2]{dr2} V. Drinfeld, On Almost Cocommutative Hopf Algebras,
{\em
Leningrad Mathematics Journal} {\bf 1} (1990), 321-342.
\bibitem
[Dr3]{dr3} V. Drinfeld, Constant quasiclassical solutions
of the quantum
Yang-Baxter equation, {\em Dokl. Acad. Nauk. SSSR} {\bf
273} {\bf No. 3}
(1983), 531-535.
\bibitem
[EG1]{eg1} P. Etingof and S. Gelaki, Some Properties of
Finite-Dimensional
Semisimple Hopf Algebras, {\em Mathematical Research
Letters} {\bf 5}
(1998), 191-197.
\bibitem
[EG2]{eg2} P. Etingof and S. Gelaki, A Method of
Construction of
Finite-Dimensional Triangular Semisimple Hopf Algebras, {\em
Mathematical
Research Letters} {\bf 5} (1998), 551-561.
\bibitem
[EG3]{eg3} P. Etingof and S. Gelaki, The Representation
Theory of
Cotriangular Semisimple Hopf Algebras, {\em International
Mathematics
Research Notices} {\bf 7} (1999), 387-394.
\bibitem
[EG4]{eg4} P. Etingof and S. Gelaki, The Classification
of Triangular  
Semisimple and Cosemisimple Hopf Algebras Over an
Algebraically Closed 
Field, {\em International Mathematics Research Notices}
{\bf 5} (2000), 223-234.
\bibitem
[EG5]{eg5} P. Etingof and S. Gelaki, On Finite-Dimensional
Semisimple and
Cosemisimple Hopf Algebras In Positive Characteristic, {\em
International Mathematics Research Notices} {\bf 16}
(1998), 851-864.
\bibitem
[EG6]{eg6} P. Etingof and S. Gelaki, On Cotriangular Hopf
Algebras,
{\em submitted}, math.QA/0002128.
\bibitem
[EGGS]{eggs} P. Etingof, S. Gelaki, R. Guralnick and J.
Saxl, Biperfect Hopf
Algebras, {\em Journal of Algebra}, to appear, math.QA/9912068.
\bibitem
[ES]{es} P. Etingof and O. Schiffmann, Lectures on
Quantum
Groups, Lectures in Mathematical Physics, {\em
International Press, Boston, MA} (1998).
\bibitem
[ESS]{ess} P. Etingof, T. Schedler and A. Soloviev,
Set-Theoretic
Solutions to the Quantum
Yang-Baxter Equation, {\em Duke Math. J.}, to appear,
math.QA/9801047.
\bibitem 
[G1]{g1} S. Gelaki, On Pointed Ribbon Hopf Algebras,
{\em Journal 
of Algebra} {\bf 181} (1996), 760-786.
\bibitem 
[G2]{g2} S. Gelaki, Quantum Groups of
Dimension $pq^2,$ {\em Israel Journal of Mathematics}
{\bf 102} (1997),
227-267.
\bibitem 
[G3]{g3} S. Gelaki, Pointed Hopf Algebras and Kaplansky's 10th 
Conjecture, {\em Journal of Algebra} {\bf 209} (1998), 635-657. 
\bibitem 
[G4]{g4} S. Gelaki, Some Examples and Properties of
Triangular Pointed Hopf Algebras, {\em Mathematical Research
Letters} {\bf 6} (1999), 563-572; see corrected version at 
math.QA/9907106.
\bibitem 
[G5]{g5} S. Gelaki, Semisimple Triangular Hopf Algebras
and Tannakian Categories, {\em submitted}.
\bibitem
[HI]{hi} R.B. Howlett and I.M. Isaacs, On Groups of
Central Type,
{\em Mathematische Zeitschrift} {\bf 179} (1982), 555-569.
\bibitem
[KaG]{kacg} G.I. Kac, Extensions of groups to ring groups,
{\em Math. USSR sbornik} {\bf 5} No. 3 (1968).
\bibitem
[KaV]{kav} V. Kac, Lie superalgebras, {\em Advances in Math.}
{\bf 26}, No.1, 1977.
\bibitem
[Kap]{kap} I. Kaplansky, Bialgebras, University of Chicago,
1975.
\bibitem
[Kash]{kash} Y. Kashina, Classification of Semisimple 
Hopf Algebras of Dimension $16,$ {\em J. Algebra}, to 
appear, math.QA/0004114.
\bibitem 
[Kass]{kass} C. Kassel, Quantum Groups, Springer, New
York, 1995.
\bibitem
[Ko]{ko} B. Kostant, Graded manifolds, graded Lie theory,
and prequantization, Differ. geom. Meth. math.
Phys., Proc. Symp. Bonn 1975, {\em Lect. Notes Math.}
{\bf 570} (1977), 177-306.
\bibitem
[L]{l} G. Lusztig, Finite Dimensional Hopf
Algebras Arising from
Quantized Universal Enveloping Algebras, {\em J. of the
A.M.S} Vol.3, No.1 (1990), 257-296.
\bibitem 
[LR1]{lr1} R. Larson and D. Radford,
Semisimple 
Cosemisimple Hopf Algebras, {\em American Journal of
Mathematics} 
{\bf 110} (1988), 187-195.
\bibitem [LR2]{lr2} R.G. Larson and D.E. Radford,
Finite-Dimensional
Cosemisimple Hopf Algebras in Characteristic $0$ are
Semisimple, {\em J.
Algebra} {\bf 117} (1988), 267-289.
\bibitem
[Ma1]{ma1} S. Majid, Foundations of quantum group theory,
Cambridge University Press, 1995.
\bibitem
[Ma2]{ma2} S. Majid, Physics for algebraists:
Non-commutative and non-cocommutative Hopf algebras by a
bicrossproduct construction, {\em J. Algebra} {\bf 130}
(1990), 17-64.
\bibitem
[Mon]{mon} S. Montgomery, Hopf algebras and their actions
on rings, {\em CBMS Lecture Notes} {\bf 82}, AMS, 1993.
\bibitem
[Mov]{mov} M. Movshev, Twisting in group algebras of finite
groups, {\em
Func. Anal. Appl.} {\bf 27} (1994), 240-244.
\bibitem 
[NZ]{nz} W.D. Nichols and M.B. Zoeller, A Hopf algebra freeness
theorem, {\em American Journal of Mathematics} {\bf 111}
(1989),
381-385.
\bibitem 
[PO1]{po1} F. Panaite and F.V. Oystaeyen,
Quasitriangular
structures for some pointed Hopf algebras of dimension
$2^n,$ {\em
Communications in Algebra}, to appear.
\bibitem 
[PO2]{po2} F. Panaite and F.V. Oystaeyen, Clifford-type
algebras as
cleft extensions for some pointed Hopf algebras, {\em
Communications in
Algebra}, to appear. 
\bibitem 
[R1]{r1} D.E. Radford, The Structure of Hopf Algebras with a 
Projection, {\em Journal of Algebra} {\bf 2} (1985), 322-347.
\bibitem 
[R2]{r2} D.E. Radford, Minimal quasitriangular Hopf algebras, 
{\em Journal of Algebra} {\bf 157} (1993), 285-315.
\bibitem 
[R3]{r3} D.E. Radford, The trace function and Hopf
Algebras {\em Journal of Algebra} {\bf 163} (1994), 583-622.
\bibitem 
[R4]{r4} D.E. Radford, On Kauffman's Knot Invariants
Arising from
Finite-Dimensional Hopf algebras, Advances in Hopf
algebras, 158, 
205-266, lectures notes in pure and applied mathematics,
Marcel Dekker, 
N. Y., 1994.
\bibitem
[Se]{se} J-P. Serre, Local Fields {\em Graduate Texts in
Mathematics} {\bf 67}.
\bibitem
[Sw]{sw} M. Sweedler, Hopf Algebras, Benjamin Press,
1968. 
\bibitem
[T]{t} M. Takeuchi, Matched pairs of groups and bismash
products of Hopf algebras, {\em Comm. Algebra} {\bf 9},
No. 8 (1981), 841-882.
\bibitem 
[TW]{tw} E.J. Taft and R. L. Wilson, On antipodes in
pointed Hopf 
algebras, {\em Journal of Algebra} {\bf 29} (1974), 27-32.
\end{thebibliography}
\end{document}